\documentclass[12pt,reqno]{amsart}
\usepackage{tikz}

\usepackage{curves}

\usepackage[small,nohug,heads=littlevee]{diagrams}
\diagramstyle[labelstyle=\scriptstyle]
\usepackage{amsmath,amsthm,amscd,amsfonts,amssymb,euscript}
\usepackage{graphics}
\usepackage{pdfpages}
\usepackage{mathtools}

\usepackage[all]{xypic}
\usepackage{latexsym}
\usepackage{color} 
\theoremstyle{plain}

\font\manual=manfnt
\def\dbend{{\manual\char127}} % dangerous bend sign

\def\danger{\begin{trivlist}\begin{footnotesize}\item[]\noindent%
\begingroup\hangindent=3pc\hangafter=-2%\clubpenalty=10000%
\def\par{\endgraf\endgroup}%
\hbox to0pt{\hskip-\hangindent\dbend\hfill}\ignorespaces}
\def\enddanger{\par\end{footnotesize}\end{trivlist}}

\def\ddanger{\begin{trivlist}\begin{footnotesize}\item[]\noindent%
\begingroup\hangindent=3pc\hangafter=-2%\clubpenalty=10000%
\def\par{\endgraf\endgroup}%
\hbox to0pt{\hskip-\hangindent\dbend\kern2pt\dbend\hfill}\ignorespaces}
\def\endddanger{\par\end{footnotesize}\end{trivlist}}

\DeclareFontFamily{OT1}{rsfs}{}
\DeclareFontShape{OT1}{rsfs}{n}{it}{<-> rsfs10}{}
\DeclareMathAlphabet{\mathscr}{OT1}{rsfs}{n}{it}
\DeclareMathOperator{\Spec}{Spec}

\newcommand{\eD}{{\mathscr D}}

\newcommand{\eE}{{\mathscr E}}
\newcommand{\eG}{{\mathscr G}}
\newcommand{\eH}{{\mathscr H}}

\newcommand{\cA}{{\mathcal A}}
\newcommand{\eV}{{\mathscr V}}

\newcommand{\eC}{{\mathscr C}}
\newcommand{\eL}{{\mathscr L}}
\newcommand{\eM}{{\mathscr M}}
\newcommand{\eN}{{\mathscr N}}

\begin{document}
\input{amssym.def}
%%%%%%%%%% Start TeXmacs macros
%\newcommand{\tmstrong}[1]{\textbf{#1}}
%%%%%%%%%% End TeXmacs macros

\numberwithin{equation}{section}

\newtheorem{guess}{\sc Theorem}[section]
\newcommand{\bth}{\begin{guess}$\!\!\!${\bf }~}
\newcommand{\eeth}{\end{guess}}

\newtheorem{propo}[guess]{\sc Proposition}%[section]

\newcommand{\bprop}{\begin{propo}$\!\!\!${\bf }~}
\newcommand{\eprop}{\end{propo}}

\newtheorem{lema}[guess]{\sc Lemma}%[section]
\newcommand{\blem}{\begin{lema}$\!\!\!${\bf }~}
\newcommand{\elem}{\end{lema}}

\newtheorem{defe}[guess]{\sc Definition}%[section]
\newcommand{\bdefe}{\begin{defe}$\!\!\!${\it }~}
\newcommand{\edefe}{\end{defe}}

\newtheorem{coro}[guess]{\sc Corollary}%[section]
\newcommand{\bcor}{\begin{coro}$\!\!\!${\bf }~}
\newcommand{\ecor}{\end{coro}}

\newtheorem{rema}[guess]{\it Remark}%[section]
\newcommand{\brem}{\begin{rema}$\!\!\!${\it }~\rm}
\newcommand{\erem}{\end{rema}}

\theoremstyle{remark}
\newtheorem{assump}[guess]{Assumption}

\newcommand{\spec}{{\rm Spec}\,}
\newtheorem{notation}{Notation}[section]
\newcommand{\bnot}{\begin{notation}$\!\!\!${\bf }~~\rm}
\newcommand{\enot}{\end{notation}}

\newcommand{\bpr}{\begin{proof}}
\newcommand{\epr}{\end{proof}}

\numberwithin{equation}{subsection} 
\newcommand{\beqa}{\begin{eqnarray}}
\newcommand{\eeqa}{\end{eqnarray}}
\newtheorem{thm}{Theorem}
\theoremstyle{definition}
\newtheorem{say}[guess]{\bf}
\newtheorem{hint}[thm]{Hint}
\newtheorem{example}[thm]{Example}
\newcommand{\bsem}{\begin{say}$\!\!\!${\it }~~\rm}
\newcommand{\esem}{\end{say}}
\newtheorem{observe}[subsubsection]{Observation}
\newsymbol \bulletledarrowleft 1309

\newcommand{\ha}{\sf h}
\newcommand{\g}{\sf g}
\newcommand{\ta}{\sf t}
\newcommand{\s}{\sf s}
\newcommand{\ctext}[1]{\makebox(0,0){#1}}
\setlength{\unitlength}{0.1mm}

\newcommand{\wt}{\widetilde}
\newcommand{\Lr}{\Longrightarrow}
\newcommand{\Aut}{\mbox{{\rm Aut}$\,$}}
\newcommand{\ul}{\underline}
\newcommand{\ol}{\bar}
\newcommand{\lr}{\longrightarrow}
\newcommand{\sh}{{\sf h}}

\newcommand{\ba}{{\mathbb A}}

\newcommand{\bc}{{\mathbb C}}
\newcommand{\bp}{{\mathbb P}}
\newcommand{\bz}{{\mathbb Z}}
\newcommand{\bq}{{\mathbb Q}}
\newcommand{\bn}{{\mathbb N}}
\newcommand{\bg}{{\mathbb G}}
\newcommand{\br}{{\mathbb R}}
\newcommand{{\bh}}{{\mathbb H}}
\newcommand{\bo}{{\bar \omega}'}
\newcommand{\po}{{\omega}'}

\newcommand{\ct}{{\mathcal T}}
\newcommand{\cc}{{\mathcal C}}
\newcommand{\cl}{{\mathcal L}}
\newcommand{\cv}{{\mathcal V}}
\newcommand{\cf}{{\mathcal F}}
\newcommand{\cb}{{\Lambda}}
\newcommand{\ch}{{\mathcal H}}

\newcommand{\mfc}{{\sf C}}
\newcommand{\ce}{{\mathcal E}}
\newcommand{\co}{{\mathcal O}}

\newcommand{\cm}{{\mathcal M}}

\newcommand{\cs}{{\mathcal S}}
\newcommand{\cg}{{\mathcal G}}
\newcommand{\ca}{{\mathcal A}}
\newcommand{\hra}{\hookrightarrow}
\newcommand{\mfu}{{\sf U}}

\newtheorem{ack}{\it Acknowledgments}       
\renewcommand{\theack}{} 

\title{A degeneration of moduli of Hitchin pairs}

\author{V. Balaji} 
\address{Chennai Mathematical Institute SIPCOT IT Park, Siruseri-603103, India,
balaji@cmi.ac.in}
%\and
\author{P. Barik} 
\address{Chennai Mathematical Institute SIPCOT IT Park, Siruseri-603103, India,
pabitra@cmi.ac.in}
\author{D.S. Nagaraj} 
\address{Institute of Mathematical Sciences, Taramani, Chennai-600115, India, dsn@imsc.res.in}

%\subjclass[2002]{Preliminary Version}
%\date{March 25, 2013}
\dedicatory{To C.S. Seshadri with respect and admiration}
\thanks{The research of the first author was partially supported by the J.C. Bose Fellowship.}
\maketitle
\small
\tableofcontents
\normalsize

\vspace{2mm}
\noindent

\section{Introduction}
The geometry of Hitchin pairs or Higgs bundles has been extensively studied for over twenty five years. The problem of  constructing a natural theory of degenerations of the moduli space of Hitchin pairs on smooth curves is therefore of some significance. The purpose of this paper is to develop such a theory. 

Let $R$ be a discrete valuation ring with quotient field $K$ and residue field an algebraically closed field $k$. Let $S = \spec ~R$, and $\spec~K$ the generic point and let $s$ be the closed point of $S$. Let $X \to S$ be a proper, flat family with generic fibre $X_{_K}$ a smooth projective curve of genus $g \geq 2$ and with closed fibre $X_{_s}$ a stable singular curve $C$ with a single node $p \in C$. Let $(n,d)$ be a pair of integers such that $gcd(n,d) = 1$. 

Let ${\eM}(n,d)^{^{H}}_{_K}$ be the moduli space of stable Hitchin pairs of rank $n$ and degree $d$ on the generic fibre $X_{_K}$  of $X/S$ . In this paper, we construct and study a degeneration of the moduli space ${\eM}(n,d)^{^{H}}_{_K}$ of rank $n$ and degree $d$ with analytic normal crossing singularities. Central to this theory is the geometry of the Hitchin fibre which reveals a somewhat new aspect of the theory of compactifications of Picard varieties of curves, which at the same time yields a degeneration of the classical Hitchin picture. In contrast to the usual theory of Picard compactifications, the ones which arise here have analytic normal crossing singularities; recall that when  the number of nodes of the curve is strictly bigger than $1$, the singularities of the compactified Picard variety is a {\em product of normal crossing singularities} and therefore {\em not} a normal crossing singularities (cf. \cite[Page 595]{caporaso}, \cite[Page 262, I]{ictp}). In this process a very natural toric picture shows up, which in a certain sense underlies the so-called {\em abelianization} philosophy (Theorem \ref{three}).

Let $C$ be a  projective curve of genus $g \geq 2$ over  $k$ and let $\cl$ be a line bundle on $C$. A Hitchin pair $(E,\theta)$, comprises of a torsion-free $\co_{_C}$-module $E$ together with a $\co_{_C}$-morphism $\theta:E \to E \otimes \cl$ called the Higgs structure. The case when $C$ is smooth and when  $\cl$ is the dualizing sheaf $\omega_{_C}$ was studied by Nigel Hitchin in the classic paper \cite{hitchin}. Hitchin gave an analytic construction of the moduli space and showed the properness of the Hitchin map; he also observed that the notion of spectral curves  comes up naturally in the theory of Hitchin pairs and gives an {\em abelianization} of the non-abelian moduli space of stable vector bundles of rank $2$ and odd degree. This theory on smooth projective curves has been generalized on numerous fronts, the most significant one being the work of Simpson (\cite{sim1}, \cite{sim2} and \cite{sim}) where it is carried through for higher dimensional smooth projective varieties. The paper by Nitsure \cite{nitsure} was the first one to give a purely algebraic construction of the moduli space and he considered the more general situation of taking an arbitrary line bundle for the Higgs structure instead of the dualizing sheaf.

 A natural approach to construct a degeneration would be to consider torsion-free Hitchin pairs on the family $X/S$ which is what is done when there is no Higgs structure involved. Quite surprisingly, unlike the usual case (where there is no Higgs structure), this method does not quite give a degeneration with the essential properties that one desires, namely flatness and even more importantly, an analogue of the Hitchin map (which is only a rational map here). We therefore develop a theory  analogous to that of the Gieseker construction of Hilbert stable bundles on a degenerating family (cf. \cite{gies}) which yields a flat degeneration (Theorem \ref{one}) and at the same time resolves the rational Hitchin map defined on the torsion-free moduli (Theorem \ref{two}). 
 
Recall that Gieseker's construction (when rank is $2$) was based on the concept of $m$-Hilbert stability, which is a GIT stability condition for points in a certain Hilbert scheme of embeddings of curves in a Grassmannian. This is not quite amenable when one goes to ranks bigger than $2$ and more so in the setting of Hitchin pairs. Nonetheless, an approach along the lines of \cite{ns2} and \cite{schmitt} does work. The process becomes somewhat intricate and the analogy needs to be delicately carried hand in hand with the moduli of relative torsion-free Hitchin pairs on $X/S$.  

The main sources for our tools are the papers by Simpson (\cite{sim1}, \cite{sim2}), and Nagaraj-Seshadri \cite{ns2} (see also Schmitt \cite{schmitt}).

Let $X \to S$ be a proper and flat fibered surface over $S = \spec (R)$, where $R$ is a  local ring of a smooth curve over $k$ with generic fibre a smooth projective curve of genus $g$ and closed fibre a singular curve $C$ with a single node $p \in C$; assume that $X$ is regular over $k$. Let $\cl$ be a relative  line bundle on $X$ and we assume that $deg(\cl|_{_C}) > deg(\omega_{_C})$, where $\omega_{_C}$ is the dualizing sheaf on $C$. The assumptions on $\cl$ are essential only for the flatness of the degeneration. For much of the existence results, we do not need any ampleness assumptions on $\cl$.
Our principal results are the following:
\bth\label{one}  
\begin{enumerate}
\item There is a quasi-projective \text{$S$-scheme} ${\eG}_{_{S}}^{^{H}}(n,d)$ of Gieseker-Hitchin pairs which is flat over $S$ and regular over $k$, with the closed fibre a divisor with (analytic) normal crossing singularities. 
\item The generic fibre is isomorphic to the classical Hitchin space ${\eM}_{_{K}}^{^{H}}(n,d)$.
\end{enumerate}
\eeth

\bth\label{two}~
\begin{enumerate}
\item We have a Hitchin map ${\sf g}_{_S}:{\eG}_{_{S}}^{^{H}}(n,d) \to {\ca}_{_S}$ to an affine space over $S$ which is {\em proper}. 

\item To a general section $\xi:S \to {\ca}_{_S}$ we can associate a spectral fibered surface $Y_{_\xi}$ over $S$ with smooth projective generic fibre $Y_{_{\xi,K}}$ and whose closed fibre $Y_{_{\xi,s}}$ is an irreducible vine curve with $n$-nodes. 

\item Let ${P_{_{\delta,Y_{_\xi}}}}$ denote the compactified relative Picard \text{$S$-scheme}  of the spectral fibered surface $Y_{_\xi}$ over $S$ (as constructed by Caporaso \cite{caporaso}). Then we have a proper birational morphism:
\beqa\label{nu}
\nu_{_*}:{\sf g}_{_S}^{-1}(\xi) \to {P_{_{\delta,Y_{_\xi}}}}
\eeqa
which is an isomorphism over the generic fibre and this map coincides with  the classical Hitchin isomorphism of the Hitchin fibre with the Jacobian of $Y_{_{\xi,K}}$.
\item The \text{$S$-scheme} ${\sf g}_{_S}^{-1}(\xi)$ gives a compactification of the Picard variety, whose fibre over $s$ is a divisor with {\em analytic normal crossing singularities}. 

\end{enumerate}
\eeth
\noindent
The fibres of the morphism $\nu_{_*}$ to the compactified Picard variety of the {\em vine curve} $Y_{_{\xi,s}}$ gets the following description:
\bth\label{three}  The morphism $\nu_{_*}$ is an isomorphism on the over subscheme of locally free sheaves of rank $1$ and for each $j$, over the stratum ${P_{_{\delta,Y_{_{\xi,s}}}}}(j)$ (see  \eqref{picstrata1}) the fibres are canonical toric subvarieties of the wonderful compactification $\overline{PGL(j)}$ obtained from the closures of the maximal tori of $PGL(j)$ (for details see Theorem \ref{relative abelianization}). These are toric varieties associated to the Weyl chamber of $PGL(j)$ (\cite{procesi}).

\eeth

The theory generalizes without serious difficulty to reducible curves, which is the content of the last section. In principle it should generalize to any stable curve but the details need to be worked out.

The layout of the paper is as follows: in Section 2 we quickly rework Simpson's theory for the family $X/S$. In Section 3 we introduce the Gieseker-Hitchin functor; in Section 5 and 6, the coarse moduli  space ${\eG}(n,d)^{^{H}}_{_{S}}$ of stable Gieseker-Hitchin pairs is constructed and we define the Gieseker-Hitchin map and prove its {\em properness}.  Section 7 and 8 are devoted to the study of the geometry of the general fibre of the Gieseker-Hitchin map for a singular curve $C$ and we conclude with the main theorem of the paper. The final section shows how the results can be generalized for a reducible curve with a single node.

\ack We are grateful to C.S. Seshadri for his interest in this work and the numerous comments and suggestions which have been invaluable.

\section{Hitchin pairs on nodal curves}
We will assume for the most part of this paper that the curve $C$ is an irreducible nodal curve  with one node over an algebraically closed field $k$. 

\subsubsection{Hitchin pairs on $C$} Let $E$ be a coherent $\co_{_C}$--module. Recall that $E$ has {\em depth} $1$ (at each $x \in C$) if and only if $E$ is of {\em pure dimension} $1$, i.e for all nonzero $\co_{_C}$--submodules $\cf \subset E$, $dim(Supp(\cf)) = 1$. Recall that if $C$ is a singular curve with nodal singularities, {\em a torsion-free $\co_{_C}$-module} is the same thing as a coherent $\co_{_C}$-module which is of {\em depth} $1$. Let $\cl$ be an invertible sheaf over $C$.

\bdefe A {\em Hitchin pair} $(E,\theta)$, comprises of a torsion-free $\co_{_C}$-module $E$ together with a $\co_{_C}$-morphism $\theta :E \to E \otimes \cl$. The map  $\theta$ is called a Higgs structure on $E$.\edefe

Let $\mu(E):= \frac{\chi(E)}{rk(E)}, \text{if}~ rk(E) \neq 0$. As $dim(C)= 1$ for us, we see immediately that
\beqa\label{pandmuirr}
\frac{{\tt p}(E,m)}{rk(E)} = m \cdot deg(\co_{_C}(1)) + \mu(E).
\eeqa
where ${\tt p}(E,m) := \chi(E \otimes \co_{_C}(m))$.

\bdefe\label{semistability} 
 A Hitchin pair $(E,\phi)$, with $E$ of constant rank, is called {\em $\mu$-(semi)stable}  if
$\mu(E_1) \leq \mu(E)$  (resp. $\lneq $), $\forall$ proper subsheaves $E_1 \subset E$ such that $\phi (E_1) \subset E_1 \otimes \cl $.
\edefe

\subsubsection{The moduli spaces}  
We will in fact work with a family of smooth curves degenerating to  an irreducible nodal curve $C$ with one node. More precisely, let  $S = \spec(R)$, with $R$ a discrete valuation ring (which is the local ring of a smooth curve over $k$) with quotient field $K$ and residue field $k$. We will denote by $s \in S$ the closed point  and $\zeta \in S$ the generic point. Let $f: X \to S$ be an $S$--scheme of relative dimension $1$, such that the generic fibre $X_{_\zeta}$ is a smooth projective curve of genus $g = p_{_a}(C) \geq 2$ and the closed fibre $X_{_s} = C$. We will also assume that $X$ is regular as a scheme over $k$. Fix  $\cl$ an arbitrary invertible sheaf on $X$. We will make no ampleness assumptions on $\cl$ till later, when we need it.

The aim of this section is to quickly summarize the  construction of the moduli  space ${\eM}^{^{H}}_{_{S}}(n,d)$ of stable Hitchin pairs on $X$ of rank $n$ and degree $d$, as a quasi-projective scheme over $S$. This is done following Simpson \cite{sim1} and \cite{sim2}; in fact, we need both the constructions in \cite{sim2}, the one in terms of $\Lambda$-modules as well as the one with pure sheaves.

\subsubsection{The moduli space of Hitchin pairs} Let $f: X \to S$, be as above.
Let $\Lambda = Sym(\cl^{*})$ as a sheaf of $\co_{_X}$-algebras. A {\em Hitchin pair} on $X$ over $S$ is a coherent $\co_{_X}$-module $E_{_S}$ together with a $\co_{_X}$-morphism $\theta_{_S}:E_{_S} \to E_{_S} \otimes \cl$. Giving $\theta_{_S}$ is equivalent to giving $\theta_{_S}:\cl^{*} \to {\mathcal E}nd(E_{_S})$, or again as an algebra homomorphism $\theta_{_S}:\Lambda \to {\mathcal E}nd(E_{_S})$ (cf. \cite[page 15]{sim2}). That is, $\theta_{_S}$ gives $E_{_S}$ a structure of a sheaf of modules over $\underline{\spec}(\Lambda)$ which is the total space of $\cl$ (see \cite[Lemma 2.13]{sim1}). By a $\Lambda$-module, we will always mean a {\em coherent} $\co_{_X}$-module with a structure as above. 

Conversely ({\it loc cit}), giving a $\co_{_X}$-coherent $\Lambda$-module is equivalent to giving a coherent $\co_{_X}$-module $E$ together with an $\co_{_X}$-module map $\theta_{_S}:E_{_S} \to E_{_S} \otimes \cl$.

Let $\co_{_X}(1)$ be the relative ample line bundle on $X$ over $S$, and for each point $t \in S$, let $X_t = X \times_{_S} k(t)$ denote the fibre and $\Lambda_t$ be the restriction of $\Lambda$ to $X_t$. A $\Lambda$-module $E_{_S}$ is $\tt p$--semistable (resp $\tt p$--stable) if $E_{_S}$ is flat over $S$ and if the restrictions $E_t$ of $E_{_S}$ to $X_t$ are $\tt p$--semistable (resp. $\tt p$--stable) of pure dimension $1$. Recall (\cite{sim1}) that $E_t$ is a $\tt p$--semi(stable) $\Lambda_t$-module if it is of pure dimension $1$ (equivalently, torsion-free), and if for any $\Lambda_t$-submodule $E_1 \subset E$, with $0 < rk(E_1) < rk(E)$, there exists an $N$ such that 
\beqa
\frac{{\tt p}(E_1,m)}{rk(E_1)} \leq \frac{{\tt p}(E,m)}{rk(E)}
\eeqa  
(resp. $\lneq $), for $m \geq N$.
\brem  It follows from \eqref{pandmuirr} that for a Hitchin pair on a curve $C$, the notion of $\tt p$--(semi)stability is the same as that of $\mu$--(semi)stability.\erem

On the singular fibre $X_s = C$, the notion of $\tt p$--(semi)stability  coincides with the notion of $\mu$--(semi)stability  with respect to the polarization  given by the ample line bundle $\co_{_X}(1)|_{{X_s}}$. The notion of pure dimension $1$ is precisely the torsion-freeness of the sheaves. 

Let ${\underline{\eM}}_{_{S}}^{^{H}}(n,d)$ be  the functor which associates for every \text{$S$-scheme} $T$,  the set ${\underline{\eM}}_{_{S}}^{^{H}}(n,d)(T)$ of the equivalence classes of  families of ${\tt p}$-semistable Hitchin pairs $(E,\theta)$ on $X_{_T} := X \times_{_S} T$ with  Hilbert polynomial $P$ given by $n$ and $d$, where $(E_{_T}, \theta_{_T}) \sim (E_{_T}', \theta_{_T}')$ if there exists a line bundle $L_{_T}$ on $T$ such that $E_{_T} \simeq E_{_T}' \otimes p_{_T}^{*}(L_{_T})$ which sends $\theta_{_T}$ to $\theta_{_T}' \otimes id$.

By \cite[Theorem 4.7]{sim1}, the set ${\underline{\eM}}_{_{S}}^{^{H}}(n,d)(T)$ can be viewed as the set of isomorphism classes of $\tt p$-(semi)stable $\Lambda_{_T}$-modules on $X_{_T}$ over $T$. Thus, the functor ${\underline{\eM}}_{_{S}}^{^{H}}(n,d)$ has a coarse moduli scheme which is realized as a GIT quotient and we denote by ${\eM}_{_{S}}^{^{H}}(n,d)$; by \cite{sim1}, this is quasi-projective over $S$ whose points correspond to \text{S}-equivalence classes  (in the sense of Seshadri) of $\tt p$-semistable torsion-free Hitchin pairs of rank $n$ and degree $d$ on the fibres $X_{_t}$. Furthermore, there is an open subset of isomorphism classes of $\tt p$-stable torsion-free Hitchin pairs and one has a separated open subfunctor ${\underline{\eM}}_{_{S}}^{^{H,s}}(n,d)$. In our situation, where we have assumed $gcd(n,d) =1$, one has ${\underline{\eM}}_{_{S}}^{^{H}}(n,d) = {\underline{\eM}}_{_{S}}^{^{H,s}}(n,d)$.

\subsubsection{Spectral construction}
Let $Z = {\bp}(\cl^* \oplus \co_{_X})$, be the projective completion of the total space of $\cl$ as a scheme over $S$; let $D_{_\infty} = Z - \cl$ denote the divisor at $\infty$ and $\pi:Z \to X$ the projection which extends the map $\pi:\cl \to X$.

\blem\label{keycorresp}(cf. \cite[Lemma 6.8]{sim2}) There is a functorial correspondence between the category of Hitchin pairs $(E_{_S},\theta_{_S})$ on $X$ and the category of coherent $\co_{_Z}$-modules $\eE$ such that $Supp(\eE) \cap D_{_\infty} = \emptyset$. The sheaf $E_{_S}$ is flat over $S$ if and only if $\eE$ is flat over $S$. Further, 
\begin{enumerate}
\item $E_s$ is torsion-free if and only if $\eE_s$ is pure of dimension $1$.
\item $(E_{_S},\theta_{_S})$ is $\mu$-semistable (resp. $\mu$-stable) $\leftrightarrow$ $\eE$ on $Z$ is $\tt p$--semistable (resp. $\tt p$--stable) in the sense of Gieseker-Maruyama-Simpson.
\end{enumerate}
\elem

\bpr As we have seen earlier, since the map $\pi:\cl \to X$ is affine, giving a coherent $\eE$ on the total space $\underline{\spec}(Sym(\cl^{*}))$ of $\cl$, is equivalent to giving a coherent $\co_{_X}$-module $E_{_S} = \pi_{_*}(\eE)$ together with an action of the $\co_{_X}$-Algebra $\Lambda$. Observe that $\Lambda = \pi_{_*}(\co_{_\cl})$. By what we have observed earlier, this is equivalent to giving a Hitchin pair $(E_{_S},\theta_{_S})$ on $X$. To $E$ one associates a coherent $\eE$ on $\cl$ by letting (cf. \cite[page 362]{ega}):
\beqa
\eE = \pi^{-1}(E_{_S}) \otimes_{_{\pi^{-1}(\Lambda)}} \co_{_{\cl}}
\eeqa  
The support of $\eE$ is proper over $X$ and this condition is equivalent to saying that $\eE$ is coherent on $\cl$ and the closure of the support of $\eE$ on $Z$ does not meet $D$. If $E_{_S}$ is torsion-free, then $\eE$ is pure of dimension $1$ and conversely.

By the equivalence of categories it follows that the subobjects also correspond to each other naturally and the equivalence of the semistable objects follows as in \cite{sim2}.\epr

Choose $k$ such that $\co_{_Z}(1) = \pi^*(\co_{_{X}}(k)) \otimes \co_{_Z}(D_{_\infty})$ is ample on $Z$; therefore $\co_{_Z}(1)|_{\cl} =  \pi^*(\co_{_{X}}(k))$. This way, for any coherent $\eE$ on $Z$ whose support does not meet the divisor $D_{_\infty}$,  the Hilbert polynomials of the $\co_{_Z}$-module $\eE$ and that of $\pi_{_*}(\eE)$ differ by a scaling factor, i.e ${\tt p}(\eE, m) = {\tt p}(\pi_{*}(\eE), km)$; hence, the notions of semistability remain intact. 

Fix a polynomial ${\tt p}$ as above of degree $1$ (the relative dimension in our case) and let ${\tt p}_{_k}(m) = {\tt p}(km)$. Then by \cite[Theorem 1.19]{sim1}, we have a coarse moduli scheme ${M}(\co_{_Z}, {\tt p}_{_k})$ of $\tt p$--semistable $\Lambda$-modules on the projective $S$--scheme $Z$ with respect to the ample line bundle $\co_{_Z}(1)$ and Hilbert polynomial ${\tt p}_{_k}$. Simpson shows that  ${M}(\co_{_Z}, {\tt p}_{_k})$  is a {\em projective $S$--scheme}. Furthermore, we have an open inclusion ${\eM}_{_{S}}^{^H}(n,d) \subset {M}(\co_{_Z}, {\tt p}_{_k})$ since ${\eM}_{_{S}}^{^H}(n,d)$ parametrizes $\tt p$--semistable $\Lambda$--modules whose support avoids the divisor $D_{_\infty}$.

\section{The Gieseker-Hitchin functor} From now on ${\it gcd(n,d) = 1}$ and unless otherwise mentioned, from now on $C$ will be an irreducible nodal curve with a single node $p \in C$. 

Let $\tilde{C}$ be its normalization and let $\nu:\tilde{C} \to C$ be the normalization map and let $\nu^{-1}(p) = \{p_1,p_2\}$.
\bdefe\label{Chains} A scheme $R^{^{(m)}}$ is called a chain of projective lines if $R^{^{(m)}} = \cup_{i=1}^{m} R_i$, with $R_i \simeq \bp^1$, and if $i \neq j$, 
\beqa
R_i \cap R_j  =  \begin{cases}
\text{singleton}&   \text {if $\mid i-j\mid = 1$}\\
\emptyset& \text {otherwise}
\end{cases}
\eeqa
\edefe
\bdefe\label{standard} Let $E$ be a vector bundle  of rank $n$ on a chain $R^{^{(m)}}$. Let $E|_{_{R_i}} = \oplus_{j=1}^{n} \co(a_{ij})$. Say that $E$ is {\em standard}  if $0 \leq a_{ij} \leq 1, \forall i,j$. Say that $E$ is {\em strictly standard} if moreover, for every $i$ there is an index $j$ such that $a_{ij} = 1$.\edefe

\bdefe\label{Gieseker curve} Let $C^{^{(m)}}$ denote the semi-stable curve which is semistably equivalent to $C$, i.e $\tilde{C}$ is a component of $C^{^{(m)}}$  and if $\nu:C^{^{(m)}} \to C$ is the canonical morphism, the fibre $\nu^{-1}(p)$ is a chain $R^{^{(m)}}$ of projective lines  of length $m$ cutting $\tilde{C}$ in $p_1$ and $p_2$ (Figure \ref{Giesekercurves}).\edefe
\begin{figure}\label{Giesekercurves}
\begin{center}
\includegraphics[scale=0.40]{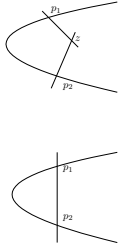} 
\caption{{\bf Semistable curve $C^{(2)}$ and $C^{(1)}$ over $C$}}
\label{Giesekercurves}
\end{center}
\end{figure} 
   
Let $p:X \to S$ be as before a family of smooth curves degenerating to the singular curve $C$. 

\bdefe\label{Gieseker surface}(cf. \cite[Definition 3.8]{kausz}) For every \text{$S$-scheme} $T$, a {\em modification} is a diagram: 
\beqa\label{modification}
\xymatrix{
X_{_T}^{^{(\text{mod})}} \ar[dr]_{p_{_T}} \ar[rr]^{\nu}& &X_{_T}\ar[dl]^{p} \\
& T  &
}
\eeqa
\begin{enumerate}
\item $p_{_T}:X_{_T}^{^{(\text{mod})}} \to T$ is {\em flat},
\item the $T$-morphism $\nu$ is finitely presented which is an isomorphism when $(X_{_T})_{_t}$ is {\em smooth}, 
\item over each closed point $t \in T$ over $s \in S$, we have $(X_{_T}^{^{(\text{mod})}}){_{_t}} = C^{^{(m)}}$ for some $m$ and $\nu$ restricts to the morphism which contracts the $\bp^1$'s on $C^{^{(m)}}$,
\end{enumerate}
\edefe
\brem We will reserve the notation $\nu$ for the modification morphism $\nu:X_{_T}^{^{(\text{mod})}} \to X_{_T}$ for any \text{$S$-scheme} $T$ and will not carry the subscript $T$ to avoid cumbersome notation. \erem
 
\bdefe\label{Giesekervectorbundle} A vector bundle $V$ on $C^{^{(m)}}$ of rank $n$ is called a {\em Gieseker vector bundle}, 
\begin{enumerate}
\item If $m = 0$, i.e $C^{^{(0)}} = C$ it is a vector bundle, else
\item if, $m \geq 1$, $V|{_{_{R^{^{(m)}}}}}$ is {\em strictly standard} and  
the direct image $\nu_{_*}(V)$ is a torsion-free on $\co_{_C}$-module. 
\end{enumerate}
A Gieseker vector bundle on a modification $X_{_T}^{(\text{mod})}$  is a vector bundle such that its restriction to each $C^{^{(m)}}$ in it is a Gieseker vector bundle. \edefe 
Let $\cl_{_{\text{mod}}}$ be the line bundle on $X_{_T}^{(\text{mod})}$ defined by $\cl_{_{\text{mod}}} := \nu^{*}(\cl)$. In particular, $\cl_{_{\text{mod}}}|{_{_{R^{^{(m)}}}}} = \co{_{_{R^{^{(m)}}}}}$, on the chain $R^{^{(m)}}$ in $C^{^{(m)}}$.

\brem\label{mayer-vietoris} If $E$ is a vector bundle on $C^{^{(m)}}$, then one has a {\it local} Mayer-Vietoris type computation to yield
\beqa
H^i(E|_{_U}) = H^i(E|_{_{V'}}) \oplus H^i(E|_{_{R^{^{(m)}}}}), \forall i \geq 1
\eeqa
where $V$ is an affine neighbourhood of the node $p \in C$ and $U = \nu^{^{-1}}(V)$ and $V' = U \cap {\tilde C}$ (cf.\cite[page 170]{ns2}).
\erem
\bdefe\label{Gieseker vb} A Gieseker-Hitchin pair on $X_{_T}^{^{(\text{mod})}}$ is a locally free Hitchin pair  $(V_{_T},\phi_{_T})$, with  an element $\phi_{_T}~\in~H^{^{0}}(X_{_T}, (p_{_T})_{_*}(\cl_{_{\text{mod}}} \otimes {\eE}nd(V_{_T}))$, i.e a $\co_{X_{_T}}$-morphism $\phi_{_T}:V_{_T} \to V_{_T}\otimes~\cl_{_{\text{mod}}}$ satisfying the following:
\begin{enumerate}
\item $V_{_T}$ is a Gieseker vector bundle on $X_{_T}^{^{(\text{mod})}}$ (Definition \ref{Giesekervectorbundle}).
\item For each closed point $t \in T$ over $s \in S$, the direct image $\nu_{_*}(V_{_t},\phi_{_t})$ is a torsion-free Hitchin pair on $X_{_t} = C$.
\end{enumerate}
A Gieseker-Hitchin pair $(V_{_T},\phi_{_T})$ is called {\em stable} if the direct image $(\nu_{_{}})_{_*}(V_{_T},\phi_{_T})$ is a family of stable Hitchin pairs on $X_{_T}$ over $T$.  
\edefe

\bdefe\label{equivalence} Two families $(V_{_T}, \phi_{_T})$ and $(V'_{_T}, \phi'_{_T})$ parametrized by $T$ are called {\em equivalent} if there exists a $X_{_T}$-automorphism $\sigma$, i.e:
\beqa\label{automor1}
\xymatrix{
X_{_T}^{^{(\text{mod})}} \ar[dr]_{\nu} \ar[rr]^{\sigma}& &X_{_T}^{^{(\text{mod})}}\ar[dl]^{\nu} \\
& X_{_T}  &
}
\eeqa
and a line bundle $\eD_{_T}$ on the parameter space $T$ such that 
\beqa
\sigma^{^*}((V_{_T}, \phi_{_T}) \otimes \eD_{_T} \simeq (V'_{_T}, \phi'_{_T}).
\eeqa
Equivalently, for each closed point $t \in T$ over $s \in S$, if there exists an automorphism $g$ of $C^{^{(m)}}$ which is identity on the normalization $\tilde{C}$, $g^{^*}(V_{_t},\phi_{_t}) \simeq (V'_{_t}, \phi'_{_t})$. 
\edefe 
\brem This notion of equivalence is the natural generalization of the notion defined in \cite[Definition 4, page 177]{ns2}.\erem

\bdefe The Gieseker-Hitchin functor ${\underline{\eG}}_{_{S}}^{^{H}}(n,d)(T)$ is defined as follows: for every \text{$S$-scheme} $T$,
\beqa
{\underline{\eG}}_{_{S}}^{^{H}}(n,d)(T) := \big[X_{_T}^{^{(\text{mod})}},(V_{_T},\phi_{_T})\big]
\eeqa
which are {\em equivalence} classes such that $(V_{_T},\phi_{_T})$ is a stable Gieseker-Hitchin pair  on $X_{_T}^{^{(\text{mod})}}$ 
and $\nu_{_*}(V_{_T},\phi_{_T}) \in {\underline{\eM}}_{_S}^{^{H}}(n,d)(T)$.   
\edefe

\section{Some auxiliary results}
Let $T = \spec~B$ be a \text{$S$-scheme} with $B$ a discrete valuation ring and let $L$ be the function field of $T$ which is assumed to be a finite extension of $K$. Assume that the closed point of $T$ maps to $s \in S$. Let $X_{_T} = X \times_{_S} T$ and let $p \in X_{_T}$ be the node; let $U$ be a formal neighbourhood of $p$ in $X_{_T}$. We recall (\cite[Page 191]{ns2} that $U$ is {\em normal} with an isolated singularity at $p$ of type $A$. By the generality of $A$-type singularities, one can realize $U$ as a cyclic quotient of the affine plane and we can write $U = \spec~C$ where $C$ is:
\beqa
C = \frac{k[[X_{_1},X_{_2},X_{_3}]]}{(X_{_1} X_{_2} - X_{_3}^{^r})}
\eeqa
i.e. $U$ is a fibered surface over $\spec~k[[X_{_3}]]$.

Let $E_{_T}$ be a family of torsion-free sheaves on $X_{_T}$ which is locally free on $X_{_L}$. Let us denote the restriction of $E_{_T}$ to $U$ by $F$. Then we have the following general lemma from \cite{ns2} which gives the complete description of such $F$'s. Recall that one can check that in this situation $F$ is not merely torsion-free but it is
also reflexive.
\blem\label{localdescrip} Let $F$ be a reflexive $C$-module which is free over the generic fibre. The $F$ is isomorphic as a $C$-module to a direct sum of ideal sheaves of the following kind:
\beqa
F \simeq  \co_{_U}^{^j} \oplus \bigoplus_{i = 1}^{l} \big(X_{_1}, X_{_3}^{^{\alpha_{_i}}}\big)^{^{\oplus r_{_i}}}.
\eeqa
where $\big(X_{_1}, X_{_3}^{^{\alpha_{_i}}}\big)$ are ideals generated by two elements.
\elem

With this lemma in place we have the following key proposition. 
\bprop\label{localprop} Let $E_{_T}$ be a reflexive sheaf on $X_{_T}$. Then there exists a modifications $\nu: X_{_T}^{^{(\text{mod})}} \to X_{_T}$, such that $V_{_T} = \nu^{^*}(E_{_T})/{tors}$ obtained by going modulo torsion is locally free and is a Gieseker vector bundle; moreover, $\nu_{_*}(V_{_T}) \simeq E_{_T}$. \eprop
\bpr The bundle $E_{_T}$ is locally free outside the single singularity $p \in X_{_T}$, therefore we reduce to the case where we take a local formal neighbourhood $U$ of $p$ and we need to show that if $F$ is reflexive sheaf on $U$, then there exists a modifications $\nu: U^{(mod)} \to U$, such that $V = \nu^{^*}(F)/{tors}$ is locally free and is a local Gieseker vector bundle and in fact, $\nu_{_*}(V) \simeq F$.

Without loss of generality, we assume in the local description Lemma \ref{localdescrip} that $j = 0$ and for simplicity of exposition we assume that the multiplicities $r_{_i} = 1, \forall i$ and 
\beqa
F \simeq \big(X_{_1}, X_{_3}^{^{a}}\big) \oplus \big(X_{_1}, X_{_3}^{^{b}}\big), ~~a < b \leq r
\eeqa
This case is sufficient to reflect the complexity of the general problem.

Let $I := \big(X_{_1}, X_{_3}^{^{a}}\big)$ and $J := \big(X_{_1}, X_{_3}^{^{b}}\big)$. Let $f:U(I) = Bl_{_I}(U) \to U$, the blow-up of $U$ with centre the ideal $I$. Its description is well-known but we need it fully. Because $U$ is integral, the scheme $U(I)$ is integral and is given as a closed subscheme of $U \times {\bp}^{^1}$. Note that $U(I)$ can be realized over $\spec C$ as ${\rm{Proj}}(Gr_{_I}(C))$, where $Gr_{_I}(C) = \bigoplus_{n \geq 0}I^{^n}$ (cf. \cite[page 4,5]{conrad}).

%\big(C[t_{_1}, t_{_2}]/(X_{_1} t_{_2} -X_{_3}^{^a} t_{_1})\big)$. 

The blow-up scheme $U(I)$ is therefore covered by two patches, $U(I) = U_{_1} \cup U_{_2}$, where 
$U_{_i} = \spec~A_{_{(t_{_i})}}, i = 1,2$. A simple computation shows that:
\beqa
A_{_{(t_{_1})}} = \frac{k[[X_{_1},X_{_2}, X_{_3}]][Y_{_2}]}{(X_{_1}~Y_{_2} - X_{_3}^{^a}, X_{_3}^{^{r-a}}Y_{_2} - X_{_2})}
\eeqa
where $Y_{_2} = \frac{t_{_2}}{t_{_1}}$ and 
\beqa
A_{_{(t_{_2})}} = \frac{k[[X_{_1},X_{_2},X_{_3}]][Y_{_1}]}{(X_{_1} - X_{_3}^{^a}Y_{_1},X_{_3}^{^{r-a}} - Y_{_1}X_{_2})}
\eeqa
where $Y_{_1} = \frac{t_{_1}}{t_{_2}}$.
Observe that $f^{^*}(I)/{tors}$ is the image of $f^{^*}(I)$ in $\co_{_{U(I)}}$; more generally, for any ideal $H$ in $\co_{_U}$ we will have the notation: 
\beqa
f^{^{\#}}(H):=  f^{^{*}}(H)\co_{_{U(I)}}.
\eeqa 
It is well-known that $f^{^{\#}}(I)$ is the relative ample bundle $\co_{_f}(1)$; explicitly this is given by $X_{_1}$, and $X_{_3}^{^a}$ on $U_{_1}$ and $U_{_2}$ respectively. We note that the generator $X_{_3}^{^a}$ in the blow-up is considered an element of homogeneous degree $1$.

Observe that since the open subsets $U_{_i}$ give the trivializing cover and the two branches $E_{_i} \subset U_{_i}, i = 1,2$, hence 
\beqa\label{forlateruse}
f^{^{\#}}(I)|_{_{E_{_i}}} = \co_{_{E_{_i}}}.
\eeqa
We now examine the inverse image $f^{^{\#}}(J)$ in $U(I)$. It is easily seen that:
\beqa
J' = \big (X_{_1}\big)
\eeqa
on $U_{_1}$ and , 
\beqa
J''= \big (X_{_3}^{^a}\big)\big(Y_{_1},X_{_3}^{^{b-a}}\big)
\eeqa
on $U_{_2}$. Thus $J$ is principal on $U_{_1}$ while its behaviour on the second patch $U_{_2}$ is similar to the initial situation, where for $C$ we have $C_{_{t_{_2}}}$ and the ideal $I$ gets replaced by the ideal $J''$ and we iterate. Let 
\beqa 
U(I,J) := Bl_{_{f^{^{\#}}(J)}}(U(I))
\eeqa 
and $g:U(I,J) \to U(I)$  the blow-up morphism. By the local nature of blow-ups, we can realize $U(I,J)$ by gluing the blow-ups of $J'$ and $J''$ on $U_{_1}$ and $U_{_2}$ respectively. 

Now since $J'$ is locally principal on $U_{_1}$, blowing it up gives us $U_{_1}$ again. Hence we can glue $Bl_{_{J''}}(U_{_2})$ with $U_{_1}$ to obtain $U(I,J) = Bl_{_J}(U(I))$ and the blow-up morphism $g:U(I,J) \to U(I)$ is obtained by gluing the blow-up morphism $g_{_{J''}}:Bl_{_{J''}}(U_{_2}) \to U_{_2}$ and the identity map on $U_{_1}$.

We now look closer at the picture which emerges in $U(I,J)$. Observe that $U(I,J) = U^{(2)}$ and $\nu = g \circ f:U^{(2)} \to U$. 

The special fibre of $U(I,J)$ over $k[[X_{_3}]]$ is the same as $U$ except that the singular point $p \in U$ is replaced by a scheme $R^{(2)} = R_{_1} \cup R_{_2}$ which is a union of two $\bp^1$'s. We make the following observations which are easily checked (see Figure \ref{blowup}):
\begin{figure}[htbp]
\begin{center}
\includegraphics[scale=0.80]{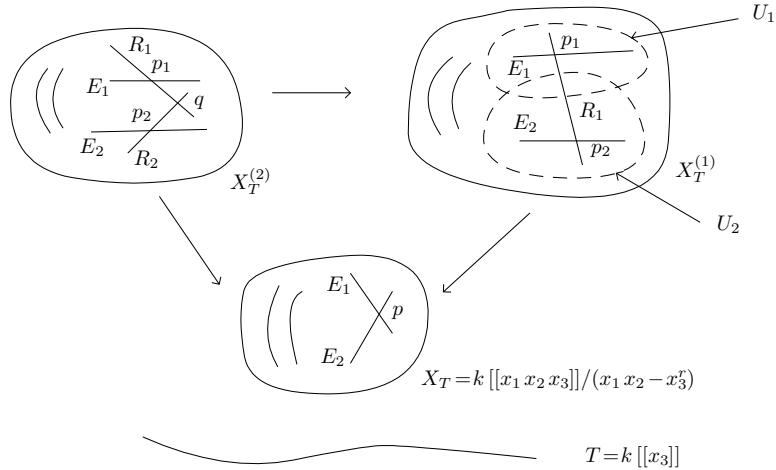} 
\caption{{\bf The blow-up picture}}
\label{blowup}
\end{center}
\end{figure}
\begin{enumerate}
\item  The {\it inverse image} $g^{^{\#}}(J)$ on $Bl_{_{J''}}(U_{_2})$ is now the relative ample $\co_{_{g_{_{J''}}}}(1)$. 
\item Hence the restriction $g^{^{\#}}(J)|_{_{R_{_{2}}}} = \co_{_{R{_{2}}}}(1)$.
\item The restriction $g^{^{\#}}(J)|_{_{R_{_{1}}}} = \co_{_{R{_{1}}}}$.
\item $g^{^*}(f^{^{\#}}(I))|_{_{R_{_{2}}}} = \co_{_{R{_{2}}}}$.
\item $g^{^*}(f^{^{\#}}(I))|_{_{R_{_{1}}}} = \co_{_{R{_{1}}}}(1)$.
\item Hence, $(g^{^{\#}} \circ f^{^{\#}})(I \oplus J)|_{_{R_{_{2}}}} = \co_{_R{_{_2}}} \oplus \co_{_R{_{_2}}}(1)$ and $(g^{^{\#}} \circ f^{^{\#}})(I \oplus J)|_{_{R_{_{1}}}} = \co_{_R{_{_1}}}(1) \oplus \co_{_R{_{_1}}}$.
\item Let $V = (g^{^{\#}} \circ f^{^{\#}})(I \oplus J)$. In other words, $V = \nu^{^{*}}(I \oplus J)/{tors}$ which is locally free on $U^{(2)}$ and on the closed fibre $R^{(2)}$ over $k[[X_{_3}]]$, it is strictly standard. 
\end{enumerate}

We need to check that the direct image $\nu_{_*}(V) = \nu_{_*}(\nu^{^{*}}(I \oplus J)/{tors})$ is firstly a flat family of torsion-free sheaves on $U$ over $k[[X_{_3}]]$ and is in fact isomorphic to $F$ we started with.

Clearly, it is enough to check this on the closed fibre of $U(I,J)$ over $k[[X_{_3}]]$. We do this recursively in the manner in which we built up the blow-up by reversing the steps. By the criterion of \cite[Lemma 2]{ns2}, we need to check that if $s$ is a section of $V|{_{_{R^{(2)}}}}$ which vanishes at $p_{_1}$ and $p_{_2}$, then it vanishes on the whole of the tree $R^{(2)}$.

As $s$ is  a section of $(g^{^{\#}} \circ f^{^{\#}})(I \oplus J)$, it can be written as $(s_{_I} \oplus s_{_J})$, such that $s_{_I}(p_{_i}) = 0 = s_{_J}(p_{_i}), ~~i=1,2$.
Now since $(g^{^{\#}} \circ f^{^{\#}})(I) = g^{^*}(f^{^{\#}}(I))$ and $f^{^{\#}}(I)$ is locally free, by the projection formula, $s_{_I}$ gives a section of $f^{^{\#}}(I)$ on $R_{_1}$ in $U(I)$ and $s_{_I}(p_{_1}) = 0$ in $R_{_1}$.

As $s_{_I}(p_{_2}) = 0$ in $R^{(2)}$, and $g^{^*}(f^{^{\#}}(I))|_{_{R_{_{2}}}}$ is $\co{_{_{R_{_2}}}}$, it vanishes everywhere on $R_{_2}$. In particular, $s_{_I}(q) = 0$, where $q = R_{_1} \cap R_{_2}$. This implies that when pushed down to $U(I)$, since $q$ maps to $p_{_2} \in R_{_1}$, the section $s_{_I}(p_{_2}) = 0 = s_{_I}(p_{_1})$, i.e $s_{_I} = 0$ since $R_{_1} \simeq \bp^{^1}$.

On the other hand, $s_{_J}(p_{_1}) = 0$ implies $s_{_J}|_{_{R_{_{1}}}} = 0$, since $g^{^{\#}}(J)|_{_{R_{_{1}}}} = \co_{_{R{_{1}}}}$; in particular, $s_{_J}(q) = 0$.  On $Bl_{_{J''}}(U_{_2})$ this is like the earlier picture, and hence $s_{_J}(q) = 0 = s_{_J}(p_{_2})$, implying $s_{_J} = 0$. This proves that $\nu_{_*}(\nu^{^{*}}(I \oplus J)/{tors})$ is firstly a flat family of torsion-free sheaves on $U$. 

Now we return to the global setting on $X_{_T}$. What we have is a Gieseker vector bundle $V_{_T}$ such that $\nu_{_*}(V_{_T})$ is a family of torsion-free sheaves on $X_{_T}$ which coincides with $E_{_T}$ on $X_{_T} - p$. Moreover, by taking double duals it follows easily that $\nu_{_*}(V_{_T})$ is reflexive (cf. \cite[page 191]{ns2}). That it is isomorphic to $E_{_T}$ now follows as in {\em loc.cit}.
\epr

\brem The existence of a Gieseker vector bundle on a modification whose direct image is the torsion-free sheaf $F$ is proven in \cite{ns2} and the result of Proposition \ref{localprop} can be deduced with a little effort from \cite{ns2}. However, the explicit form in which we prove this here is absolutely essential for the paper especially for the Hitchin pair situation. 
In fact, as we discovered later (see \cite{lipman}), if $\mathcal O$ is a $2$-dimensional analytic local rational singularity, $\nu:\tilde{\co} \to \co$ the minimal resolution of singularities, and if $M$ is a reflexive $\co$-module, then $\tilde{M} := \nu^{^*}(M)/tors$ is locally free and $\nu_{_*}(\tilde{M}) = M$. This holds in positive characteristics as well and this makes our entire set of results characteristic free.
\erem

\subsubsection{A key properness result} We make an ad hoc definition for the purposes of this paper.
\bdefe\label{horizontal properness} Let $F,G: \{\text{\text{$S$-scheme}}s\} \to \{Sets\}$ be two functors with $S = \spec~R$ as before. Suppose that $f:F \to G$ is a $S$-morphism of functors (more precisely, a natural transformation). We say, $f$ is {\em horizontally proper} if  the following valuative property holds:
let $A$ be a discrete valuation ring with function field $L$ such that $L$ is a finite extension of $K$ and $\spec~ A \to \spec~ R$ is surjective. Then for every map $\alpha \in F(L)$, if the composite $f(\alpha) \in G(L)$ extends to an element $G(A)$, then $\alpha$ also extends to an element in $F(A)$. \edefe 
This definition becomes significant along with the following observation.

\blem\label{quasiprojmorphs} Let $f:F \to G$ be a quasi-projective $S$-morphism of schemes of finite type such that $f_{_\zeta}:{F}_{_\zeta} \to {G}_{_\zeta}$ over the generic point is proper. Suppose further that the structure morphisms $F,G \to S$ are surjective and that $f$ is horizontally proper. Then $f$ is proper. \elem
\bpr The proof is essentially there in \cite[Page 188]{ns2}. Now since $f$ is quasi-projective, there is a {\em projective morphism} $\bar{f}:Z \to G$ and a diagram:
\beqa\label{}
\xymatrix{
{F} \ar[dr]_{f} \ar[rr]^{\hra}& & Z\ar[dl]^{\bar{f}} \\
& G &
}
\eeqa
As $f_{_\zeta}:{F}_{_\zeta} \to {G}_{_\zeta}$  is assumed to be proper, we may assume that $Z$ is the closure of ${F}_{_\zeta}$. To show that $Z = F$, take a point in $z \in Z_{_s}$ over the closed point $s \in S$. Choose a smooth curve $T$ connecting $z$ to a point in ${F}_{_\zeta}$, which therefore surjects onto $S$. 
Thus, we have an open subset $\Spec(L) = U \subset T$, such that $Im(U) \subset {F}$ and since $T$ maps to $Z$, it maps by composition to $G$.
 Now use the horizontal properness of $f$ to conclude that $Im(T) \subset F$ and hence $z \in F$. \epr
 
\bth\label{keyproperness} Let $\nu:{\underline{\eG}}_{_{S}}^{^{H}}(n,d) \to {\underline{\eM}}_{_S}^{^{H}}(n,d)$ be the morphism of functors induced by $\nu_{_*}$. Then the morphism $\nu$ is {\em horizontally proper}.\eeth
\bpr In the notation of the proof of Lemma \ref{quasiprojmorphs}, suppose that we have a point $g' \in {\underline{\eG}}_{_{S}}^{^{H}}(n,d)(U)$ such that $\nu(g')$ extends to a point ${\underline{\eM}}_{_S}^{^{H}}(n,d)(T)$. Then we need to show that $g'$ is the image of a point $g \in {\underline{\eG}}_{_{S}}^{^{H}}(n,d)(T)$.

Therefore, we may assume that the morphism $g'$ gives a family $(V_{_U}, \phi_{_U})$ of stable Hitchin pairs on the family of smooth curves $X_{_U}$. We need to show that there is a surface $X_{_T}^{^{(\text{mod})}}$ such that $\lim_{u \to \tau}(V_{_u}, \phi_{_u}) = (V_{_\tau}, \theta_{_\tau})$ exists as a stable Gieseker-Hitchin pair on $C^{^{(\ell)}}$ for some $\ell$.

As $\nu$ induces the identity map on $U$, we may view $(V_{_U}, \theta_{_U})$ as a family $(E_{_U}, \theta_{_U})$ on $X_{_U}$. By assumption, this extends to a torsion-free Hitchin pair $\lim_{u \to \tau}(E_{_u}, \theta_{_u}) = (E_{_\tau}, \theta_{_\tau})$ on the singular curve $C$. Let us denote this family by $(E_{_T}, \theta_{_T})$, a point of ${\underline{\eM}}_{_{S}}^{^{H}}(n,d)(T)$.

By Proposition \ref{localprop}, the family $E_{_T}$ lifts to a family $V_{_T}$ on a suitable modification $X_{_T}^{(\text{mod})}$ such that  $V_{_T}$ is a {\em Gieseker vector bundle} and also a point of ${\underline{\eG}}_{_{S}}(n,d)(T)$. The choice of the number $\ell$ which gives the length of the chain of $\bp^1$'s is dictated by the local type of the torsion-free sheaf $E_{_\tau}$. 

In fact,by Proposition \ref{localprop} we have an isomorphism
\beqa
\Big(\frac{\nu^*(E_{_T})}{tors}\Big)\simeq V_{_T}
\eeqa
Pulling back the Higgs structure $\theta_{_T}:E_{_T} \to  E_{_T} \otimes \cl$ by $\nu$, we get a Higgs structure $\nu^*(\theta_{_T}): \nu^*(E_{_T}) \to \nu^*(E_{_T}) \otimes \nu^*(\cl) = \nu^*(E_{_T}) \otimes \cl_{_\ell}$.

This gives a Higgs structure on $\Big(\frac{\nu^*(E_{_T})}{tors}\Big)$  and hence a morphism $\phi_{_T}:V_{_T} \to V_{_T}  \otimes \cl_{_\ell}$ of $\co_{X_{_T}^{(\text{mod})}}$-modules. 

To see that $\nu_{_*}(\phi_{_T}):E_{_T} \to E_{_T} \otimes \nu_{_*}(\cl_{_\ell}) = E_{_T} \otimes \cl$, gives the original Hitchin pair $(E_{_T}, \theta_{_T})$, observe that on $X_{_U}$, they give the same point of ${\underline{\eM}}_{_S}^{^{H}}(n,d)(U)$. In fact, they coincide on the whole of $X_{_T} - p$ where $E_{_T}$ is a vector bundle. As $E_{_T}$ is torsion-free, this implies immediately that they coincide on the whole of $X_{_T}$.

Clearly, $(V_{_T},\phi_{_T}))$  extends $(V_{_U},\phi_{_U})$ on $X^{^{(mod)}}_{_U}$ and furthermore, the identification $\nu_{_*}(\phi_{_T}) =  \theta_{_T}$ implies that $\phi_{_\tau}:V_{_\tau} \to V_{_\tau} \otimes (\cl_{_\ell})_{_\tau}$ give the desired limiting stable Gieseker-Hitchin pair on $C^{^{(\ell)}}$ proving the theorem.         
\epr

\section{Coarse moduli for Gieseker-Hitchin pairs} 

\subsubsection{Some deformation theory} Fix a positive integer $m$ and let $C^{^{(m)}}$ be a Gieseker curve. As before let $\nu:C^{^{(m)}} \to C$ be the canonical morphism and let ${\cl}_{_m} = \nu^{^*}(\cl)$ be the line bundle on $C^{^{(m)}}$ for defining the Higgs structure. Let $(V,\phi)$ be a Gieseker-Hitchin pair on $C^{^{(m)}}$. Our primary interest in the subsection is some remarks on the infinitesimal deformations of the pair $(V,\phi)$. 

We follow the papers \cite[page 296]{nitsure} and \cite{bisram}.
\bdefe\label{the complex} The pair $(V,\phi)$ defines a complex
\beqa\label{complex}
\eC^{^\centerdot}:= \eC^{^0} \stackrel{e(\phi)}\lr \eC^{^1} \to 0
\eeqa
where $\eC^{^0} = {\eE}nd(V)$ and $\eC^{^1} = {\eE}nd(V) \otimes {\cl}_{_m}$ and the map  
$e(\phi)$ send a local section of $s$ of ${\eE}nd(V)$ to $e(\phi)(s) = \phi \circ s - (Id_{_{\cl_{_m}}} \otimes s) \circ \phi$.
\edefe
\brem We have a short exact sequence:
\beqa
0 \to {\eE}nd(V) \otimes {\cl}_{_m}[1] \to \eC^{^.} \to {\eE}nd(V) \to 0
\eeqa
which gives the cohomology long exact sequence:\footnotesize
\beqa\label{les}
\begin{CD}
0 \to {\mathbb H}^{^0}(\eC^{^\centerdot}) \to H^{^0}({\eE}nd(V)) \to H^{^0}({\eE}nd(V)\otimes {\cl}_{_m}) \to {\mathbb H}^{^1}(\eC^{^\centerdot})\\ 
\to H^{^1}({\eE}nd(V)) \to H^{^1}({\eE}nd(V)\otimes {\cl}_{_m}) \to {\mathbb H}^{^2}(\eC^{^\centerdot}) \to 0
\end{CD}
\eeqa\normalsize
where the ${\mathbb H}^{^i}(\eC^{^\centerdot})$ are the hypercohomologies of the complex 
and the map $H^{^0}({\eE}nd(V)) \to H^{^0}({\eE}nd(V)\otimes {\cl}_{_m})$ is the map $e(\phi)$.\erem

\bprop\label{hypervanish} Let $\cl$ be a line bundle on the curve $C$ such that $deg(\cl) > deg(\omega_{_C})$. Then for any stable Gieseker-Hitchin pair $(V,\phi)$ on $C^{^{(m)}}$, with $\phi \in H^{^0}({\eE}nd(V)\otimes {\cl}_{_m})$, the hypercohomology group ${\mathbb H}^{^2}(\eC^{^\centerdot}) = 0$.\eprop
\bpr  Let $\omega_{_m}:= \omega_{_{C^{^{(m)}}}}$ denote the dualizing sheaf on $C^{^{(m)}}.$ We recall that for a Gieseker curve $C^{^{(m)}}$, one has the property 
\beqa
\omega_{_m} = \nu^{^*}(\omega_{_C})
\eeqa 
(see \cite[Lemma, page 200]{ns2}).

By the long exact sequence \eqref{les}, it suffices to show that the map 
\beqa
H^{^1}({\eE}nd(V)) \to H^{^1}({\eE}nd(V)\otimes {\cl}_{_m})
\eeqa
is surjective.  Observe that by Serre duality this map is dual to the  map:
\beqa
H^{^0}({\eE}nd(V) \otimes \omega_{_m} \otimes {\cl}_{_m}^{^{-1}}) \to H^{^0}({\eE}nd(V)\otimes \omega_{_m})
\eeqa
and hence we need to show the injectivity of this canonical map which is given as follows:
\beqa\label{inject}
s \mapsto e(\phi)(s).
\eeqa
Let $s \in H^{^0}({\eE}nd(V) \otimes \omega_{_m} \otimes {\cl}_{_m}^{^{-1}})$ i.e $s:V \to V \otimes \omega_{_m} \otimes {\cl}_{_m}^{^{-1}}$, and we suppose that $e(\phi)(s) = 0$. 

By the projection formula,  if $(\ce,\theta):= (\nu_{_*}(V),\nu_{_*}(\phi))$, we have:
\beqa
\nu_{_*}(s):\ce \to \ce \otimes \omega_{_C} \otimes \cl^{^{-1}}
\eeqa
and it is checked easily that the condition $e(\phi)(s) = 0$ translates to $e(\theta)(\nu_{_*}(s)) = 0$. By tensoring with $\omega_{_C}^{^{-1}} \otimes \cl$, the morphism $\nu_{_*}(s)$ gives rise to 
\beqa
\psi:\ce \otimes \omega_{_C}^{^{-1}} \otimes \cl \to \ce
\eeqa
and since $e(\theta)(\nu_{_*}(s)) = 0$, it follows that $Im(\psi)$ is a Higgs subsheaf of the stable Hitchin pair $(\ce,\theta)$ (since $(V,\phi)$ is assumed stable). Further, $(\ce \otimes \omega_{_C}^{^{-1}} \otimes \cl, \theta \otimes Id)$ is also a stable Hitchin pair. Comparing degrees and observing that by assumption $deg(\omega_{_C}^{^{-1}} \otimes \cl) > 0$, it follows that $\mu(Im(\psi)) > \mu(\ce)$ contradicting the stability of $(\ce,\theta)$. Thus, $Im(\psi) = 0$ and hence $\psi = 0$. In other words, we conclude that 
$\nu_{_*}(s) = 0$.

By \cite[Remark 4(i), page 176]{ns2}, the direct image $\nu_{_*}(V) = \ce$ completely determines the restriction $V|_{_{\tilde C}}$ to the normalization $\tilde C$ of $C$. Again by \cite[Remark 4(iii)]{ns2}, since $\nu_{_*}(s) = 0$, it follows that the restriction $s|{_{\tilde C}}:V|_{_{\tilde C}} \to V|_{_{\tilde C}} \otimes (\omega_{_m} \otimes {\cl}_{_m}^{^{-1}} )|_{_{\tilde C}}$ is the {\em zero map}. This implies that 
\beqa\label{onepart}
s(p_{_1}) = s(p_{_2}) = 0.
\eeqa
Viewing $s$ as a section of ${\eE}nd(V) \otimes \omega_{_m} \otimes {\cl}_{_m}^{^{-1}}$ on $C^{^{(m)}}$, and noting that $\omega_{_m} \otimes {\cl}_{_m}^{^{-1}}|{_{_{R^{^{(m)}}}}} = \co_{_{R^{^m}}}$, by restriction we get a section $s$ of ${\eE}nd(V)|_{{_{R^{^{(m)}}}}}$ such that $s(p_{_1}) = s(p_{_2}) = 0$.

If $s$ is a section of ${\eE}nd(V)|_{{_{R^{^{(m)}}}}}$ with  
$s(p_{_1}) = s(p_{_2}) = 0$ then by Lemma \ref{more2} below, it follows that $s$ vanishes.

Now as $s|{_{\tilde C}} = 0$,  we conclude that the section $s$ of 
${\eE}nd(V) \otimes \omega_{_m} \otimes {\cl}_{_m}^{^{-1}}$ obtained by gluing $s$ on $R^{^{(m)}}$ with   $s|{_{\tilde C}}$ is zero. This shows that the map \eqref{inject} is injective proving the proposition.
\epr

\blem\label{moreongiesbundles} Let $V$ be a Gieseker bundle of rank $n$ (Definition \ref{Giesekervectorbundle}) on $C^{^{(m)}}$. Let $R^{^{(m)}} = \cup_{_{j=1}}^{^m} R_{_j}$. Then 
\beqa\label{statement1}
V|_{_{R^{^{(m)}}}} = \oplus_{_{i = 1}}^{^n} \eL_{_i}
\eeqa
such that if $\eL_{_{ij}} = \eL_{_i}|{_{_{R_{_j}}}}$, then $deg(\eL_{_{ij}}) \geq 0$ and
$\sum_{j = 1}^{^{m}} deg(\eL_{_{ij}}) \leq 1$.\elem 
\bpr The first statement \eqref{statement1} follows for Gieseker vector bundles from \cite[Lemma 1]{ns2}. More generally, by \cite[Theorem 2.2]{teod}, \eqref{statement1} is true for any vector bundle on a tree. The second statement follows immediately from \cite[Proposition 5]{ns2}.\epr

\blem\label{more2} Let $V$ be a Gieseker bundle of rank $n$ on on $C^{^{(m)}}$. Then
\beqa
{\eE}nd(V) = \oplus_{_{i,k = 1}}^{^n} {\eH}om(\eL_{_i}, \eL_{_k})
\eeqa 
such that ${\eH}om(\eL_{_i}, \eL_{_k})|_{_{R_{_j}}}$ is isomorphic to $\co(1)$ for at most one index $j$ and to $\co(-1)$ for at most one index $j$.
\elem
\bpr This follows immediately from Lemma \ref{moreongiesbundles}.\epr

\brem We observe that if $\cl = \omega_{_C},$ then $dim~{\mathbb H}^{^2}(\eC^{^\centerdot}) = 1.$ \erem

\subsubsection{The total family construction}
We recall from \cite[Page 179]{ns2} the notion of  Gieseker functor with respect to a choice of a bounded family of torsion-free sheaves on $X$. In \cite{ns2} there is a slight mixing of terminologies which should become clear in our discussion. We are eventually interested in the Gieseker-Hitchin pairs therefore we begin by considering torsion-free Hitchin pairs.

We now recall from Simpson \cite[Theorem 3.8]{sim1}, the construction of a parametrizing scheme $R^{^{\Lambda,s}}_{_S}$ for stable $\Lambda$-modules with fixed Hilbert polynomial $P$. There is a choice of a positive integer $N$, such that
the functor  which   associates to each \text{$S$-scheme} $T$, the set of isomorphism classes of pairs $(\eE, \alpha)$, with $\eE$ a coherent $\Lambda$-module with Hilbert polynomial $P$ on $X_{_t}$ and for each $t \in T$:
\beqa
\alpha_{_t}:k^{^{P(N)}} \stackrel{\simeq}\longrightarrow H^{^0}(X_{_t}, \eE(N)_{_t})
\eeqa
is representable by a quasi-projective scheme $R^{^{\Lambda}}_{_S}$ over $S$. 
By the identification of torsion-free Hitchin pairs with $\Lambda$-modules, the \text{$S$-scheme} $R^{^{\Lambda}}_{_S}$ parametrizes torsion-free Hitchin pairs of rank $n$ and degree $d = P(N) + n(g-1)$ on $X/S$.  Also we have an open subset of stable points, $R^{^{\Lambda,s}}_{_S} \subset R^{^{\Lambda}}_{_S}$ and by (cf. \cite[Theorem 4.7]{sim1}),  
\beqa 
{\eM}_{_S}^{^H}(n,d) \simeq R^{^{\Lambda,s}}_{_S}/PGL(\ell)
\eeqa 
with $\ell = P(N) = dim(k^{^{P(N)}})$ and $d = \ell + n(g-1)$; and since we have assumed that $gcd(n,d) = 1$, $R^{^{\Lambda,s}}_{_S}$ is a principal $PGL(\ell)$-bundle over ${\eM}_{_S}^{^H}(n,d)$. 

We fix this total family $R^{^{\Lambda}}_{_S}$. Inside this quasi-projective variety we have a closed subscheme of torsion-free Hitchin pairs $(E,0)$ i.e. with the {\it $``0"$} Higgs structure. Let us denote this subset by $R_{_S} \subset R^{^{\Lambda}}_{_S}$. Again we have an open subset $R^{^{s}}_{_S} \subset R_{_S}$ of stable torsion-free Hitchin pairs $(E,0)$ and this is also invariant under the action of $PGL(\ell)$. Furthermore, the quotient $R^{^{s}}_{_S}/PGL(\ell) \simeq {\eM}_{_S}(n,d)$, is the moduli space of stable torsion-free sheaves on $X/S$ with rank $n$ and degree $d$ (without Higgs structure).

\subsubsection{The relative functor} We now recall from \cite[Definition 7]{ns2} the definition of the {\em Gieseker functor relative to $R_{_S}$} (in \cite{ns2}, this functor is ambiguosly called {\it Gieseker functor}!).
\bdefe\label{GiesekerfunctorreltoR} 
Let ${\eG}_{_{R_{_S}}}:(\text{\text{$S$-scheme}}s) \to (Sets)$ be the functor defined by:\small
\beqa
{\eG}_{_{R_{_S}}}(T) = \{{closed~subschemes~\Delta_{_T} \subset X \times_{_S} T\times \ Grass(\ell,n)}\}
\eeqa
\normalsize
such that:
\begin{enumerate}
\item The projection $\Delta_{_T} \to T\times \ Grass(\ell,n)$ is a closed immersion.
\item The projection $\Delta_{_T} \to T$ is a flat family of curves $\Delta_{_t}$, $t \in T$, such that $\Delta_{_T}$ is a fibered scheme over $T$ of the form $X_{_T}^{^{(\text{mod})}}$.
\item The projection $\Delta_{_T} \to X_{_T}$ is the modification map $X_{_T}^{^{(\text{mod})}} \stackrel{{\nu}} \to X_{_T}$. Furthermore, if $V_{_T}$ is the pull-back and restrictio to $\Delta_{_T}$ of the tautological quotient bundle of rank $n$  on $Grass(\ell,n)$, then $V_{_T}$ is a Gieseker vector bundle on $\Delta_{_T} = X_{_T}^{^{(\text{mod})}}$ of rank $n$ and degree $d = \ell + n(g-1)$. 
\end{enumerate}
\edefe
\brem 
We see that the direct image $\nu_{_*}(V_{_T}) = E_{_T}$ is a point of $R_{_S}(T)$.
\erem
It is shown in \cite{gies} and \cite[Proposition 8]{ns2} that the functor ${\eG}_{_{R_{_S}}}$ is represented by a $PGL(\ell)$-invariant open subscheme $\mathcal Y$ of the \text{$S$-scheme} $Hilb^{^{P_{_1}}}(X \times_{_S} Grass(\ell,n))$. Let $\Delta_{_\mathcal Y} \subset X \times_{_S} \mathcal Y\times_{_S} \ Grass(\ell,n)$ be the universal object defining the functor ${\eG}_{_{R_{_S}}}$. Then we have a canonical projection 
\beqa\label{key0}
\vartheta:\mathcal Y \to R_{_S}
\eeqa 
which for each \text{$S$-scheme} $T$, sends a point $(\Delta_{_T}, V_{_T}) = (X_{_T}^{^{(\text{mod})}}, V_{_T})$ to the direct image $\nu_{_*}(V_{_T}) = E_{_T}$, with $E_{_T} \in R_{_S}(T)$.

The imbedding $\Delta_{_\mathcal Y} \subset X \times_{_S} {\mathcal Y}\times_{_S} \ Grass(\ell,n)$  gives the natural projections of schemes over $S$:
\beqa\label{lemma4diag}
\xymatrix{
{\Delta_{_\mathcal Y}} \ar[dr]_{p} \ar[rr]^{q}& & {\mathcal Y}\ar[dl]^{r} \\
& S &
}
\eeqa
 Let $f:\Delta_{_\mathcal Y} \to X$ be the projection. Let $\cl_{_\Delta} = f^*(\cl)$ and  $\eV$ be the universal vector bundle on $\Delta_{_\mathcal Y}$ (obtained by pulling back the tautological quotient bundle from $Grass(\ell,n)$ and restricted to $\Delta_{_\mathcal Y}$). Let $\ce:= {\eE}nd(\eV) \otimes \cl_{_\Delta}$.

\bdefe\label{yh} Let $T$ be a $\mathcal Y$-scheme and $f:T \to {\mathcal Y}$. Let $q_{_T}:\Delta_{_T} \to T$ be the induced projection and $\cl_{_{\Delta_{_T}}}$ the corresponding line bundle. Let $\ce_{_T} = 
{\eE}nd(\eV) \otimes \cl_{_{\Delta_{_T}}}$. Consider the functor :
\beqa
{\eG}^{^{H}}_{_{R_{_S}}}:{\mathcal Y}-schemes \to Groups
\eeqa
which associates to each point $f \in {\mathcal Y}(T)$ , the group 
$H^{^0}(T, (q_{_T})_{_*}(\ce))$. 
\edefe

As $\mathcal Y$ is a reduced scheme (\cite{ns2}), by \cite[Lemma 3.5]{nitsure}, the functor ${\eG}^{^{H}}_{_{R_{_S}}}$ is representable; moreover, there exists a {\em linear} $\mathcal Y$-scheme ${\mathcal Y}^{^H}$  which represents it.

\brem \label{key1} By Definition \ref{yh},  for an \text{$S$-scheme} $T$, a point in ${\eG}^{^{H}}_{_{R_{_S}}}(T)$ is given by $(\Delta_{_T},V_{_T},\phi_{_T})$, where 
\begin{enumerate}
\item $(\Delta_{_T},V_{_T}) = (X_{_T}^{^{(\text{mod})}}, V_{_T}) \in {\eG}_{_{R_{_S}}}(T)$, and  
\item $(V_{_T},\phi_{_T})$ is a Gieseker-Hitchin pair on $X^{(m)}_{_T}$. 
\end{enumerate}
Thus by taking direct images and by the universal property of $R^{^{\Lambda}}_{_S}$, we see that the direct image torsion-free Hitchin pair $(\nu_{_*}(V_{_T},\phi_{_T})) = (E_{_T}, \theta_{_T})$ gives a point in $R^{^{\Lambda}}_{_S}(T)$ and we get   a canonical morphism (see \eqref{key0}) 
\beqa\label{key2}
\vartheta:{\mathcal Y}^{^H} \to R^{^{\Lambda}}_{_S}.
\eeqa 
\erem
\bprop\label{giestotfproper} The morphism $\vartheta:{\mathcal Y}^{^H} \to R^{^{\Lambda}}_{_S}$ obtained in \eqref{key2} is {\em proper}. \eprop
\bpr The schemes involved are quasi-projective schemes over $S$ and hence the morphism $\vartheta$ is a quasi-projective morphism. Also the morphism $\vartheta:{\mathcal Y}^{^H} \to R^{^{\Lambda}}_{_S}$ is  an isomorphism over the generic point of $S$. Hence by Lemma \ref{quasiprojmorphs} we need to check only {\it horizontal properness} of $\vartheta$. This follows now from Theorem \ref{keyproperness}. \epr

Let 
\beqa
{\mathcal Y}_{_{st}}^{^H}:= \vartheta^{-1}(R^{^{\Lambda,s}}_{_S}).
\eeqa
and let ${\eG}^{^{H,st}}_{_{R_{_S}}}$ be the corresponding subfunctor of ${\eG}^{^{H}}_{_{R_{_S}}}$.

Let ${\eG}'_{_{R_{_S}}}$ be the functor obtained from ${\eG}_{_{R_{_S}}}$, by forgetting the imbeddings into the Grassmannians (see \cite[Appendix, page 197]{ns2}). Then we have a 
canonical morphism
\beqa\label{formalsmoothmap}
{\eG}^{^{H,st}}_{_{R_{_S}}} \to {\eG}'_{_{R_{_S}}}
\eeqa

\bprop\label{smooth} Let $\cl$ be a line bundle on $X$ with the assumptions of Proposition \ref{hypervanish}. Then the morphism \eqref{formalsmoothmap} is formally smooth. In particular, the \text{$S$-scheme} ${\mathcal Y}_{_{st}}^{^H}$ is regular over $k$ with a divisor $({\mathcal Y}_{_{st}}^{^H})_{_s}$ with normal crossing singularities. \eprop
\bpr Let $T$ be the spectrum of an Artin local ring, and $T_{_o} \subset T$ the subscheme defined by an ideal of dimension $1$. Let $\theta \in {\eG}'_{_{R_{_S}}}(T)$ be such that the restriction $\theta_{_o} \in {\eG}'_{_{R_{_S}}}(T_{_o})$ can be lifted to an element of  
${\eG}^{^{H,st}}_{_{R_{_S}}}(T_{_o})$, then we need to show that $\theta$ itself can be lifted to an element of ${\eG}^{^{H,st}}_{_{R_{_S}}}(T)$. If $\theta$ is defined by $\Delta_{_T} \to T$, with $\Delta_{_T} \subset X \times_{_S} T \times Grass(\ell,n)$, then the lift of $\theta_{_o}$ defines a stable Gieseker-Hitchin pair $(V_{_{T_{_o}}},\phi_{_{T_{_o}}})$ on the restriction $\Delta_{_{T_{_o}}}$ of $\Delta_{_T}$ to $T_{_o}$.

The problem is to extend the pair $(V_{_{T_{_o}}},\phi_{_{T_{_o}}})$ to a stable Gieseker-Hitchin pair $(V_{_T},\phi_{_{T}})$ on $\Delta_{_T}$ as well as the sections of $V_{_{T_{_o}}}$ to those of $V_{_T}$. The second issue is taken care of as in \cite{ns2}. The key issue for us is the first one. Let $(V,\phi)$ be the restriction of $(V_{_{T_{_o}}},\phi_{_{T_{_o}}})$ to the closed fibre of $\Delta_{_{T_{_o}}}\to T_{_o}$. Then by \cite[Theorem 3.1]{bisram}, the obstruction to extending the pair $(V_{_{T_{_o}}},\phi_{_{T_{_o}}})$ lies in the second hypercohomology ${\mathbb H}^{^2}(\eC^{^\centerdot})$ of the complex $\eC^{^\centerdot}$ defined in \eqref{complex}. By Proposition \ref{hypervanish}, together with the assumptions on the line bundle $\cl$, it follows that ${\mathbb H}^{^2}(\eC^{^\centerdot}) = 0$. This implies the formal smoothness of the morphism \eqref{formalsmoothmap}. Now following the arguments in \cite[Appendix]{ns2} it follows in much the same manner that the \text{$S$-scheme} ${\mathcal Y}_{_{st}}^{^H}$ which represents the functor ${\eG}^{^{H,st}}_{_{R_{_S}}}$ is regular over $k$ and the closed fibre $({\mathcal Y}_{_{st}}^{^H})_{_s}$ has normal crossing singularities.

\epr

\bprop\label{coarsemoduli} The action of $PGL(\ell)$ on $R^{^{\Lambda,s}}_{_S}$  lifts to ${\mathcal Y}_{_{st}}^{^H}$ and the geometric quotient ${\mathcal Y}_{_{st}}^{^H} \to {\mathcal Y}_{_{st}}^{^H}/PGL(\ell)$ exists. The quotient 
\beqa
{\eG}_{_{S}}^{^{H}}(n,d):= {\mathcal Y}_{_{st}}^{^H}/PGL(\ell);
\eeqa
gives the moduli scheme for the Gieseker-Hitchin functor 
${\underline{\eG}}_{_{S}}^{^{H}}(n,d)$. \eprop
\bpr It is seen easily that  $\vartheta:{\mathcal Y}^{^H} \to R^{^{\Lambda}}_{_S}$ is generically an isomorphism and on the closed fibre it is {\em birational} over the open subset of locally free Hitchin pairs on $C$. Under these conditions, we have a natural algorithm formulated in \cite[page 179-180]{ns2} to show the existence of a quotient which is quasi-projective over $S$.

Now since $R^{^{\Lambda,s}}_{_S} \to R^{^{\Lambda,s}}_{_S}/PGL(\ell)$ is a principal $PGL(\ell)$-bundle, it follows that the GIT quotient ${\mathcal Y}_{_{st}}^{^H} \to {\mathcal Y}_{_{st}}^{^H}/PGL(\ell)$ is also a principal bundle and gives a natural quasi-projective \text{$S$-scheme} structure on ${\mathcal Y}_{_{st}}^{^H}/PGL(\ell)$. \epr
\bcor\label{giestotfproper1}
The map $\vartheta$ descends to give a {\em proper and birational} morphism
\beqa
\nu_{_*}:{\eG}_{_{S}}^{^{H}}(n,d) \to {\eM}_{_{S}}^{^{H}}(n,d).
\eeqa
\ecor
\bpr Properness follows immediately from Proposition \ref{giestotfproper}. Birationality follows from the isomorphism over smooth curves. \epr

\bcor\label{asnc} The Gieseker-Hitchin moduli space ${\eG}_{_{S}}^{^{H}}(n,d)$ of stable Gieseker-Hitchin pairs is  flat over $S$. Furthermore, as a scheme over $k$ it is regular and the closed fibre ${\eG}_{_{S}}^{^{H}}(n,d)_{_s}$ is a divisor with analytic normal crossing singularities. \ecor
\bpr  As ${\mathcal Y}_{_{st}}^{^H} \to {\mathcal Y}_{_{st}}^{^H}/PGL(\ell)$ is a principal bundle,  by Proposition \ref{smooth} it follows that
${\eG}_{_{S}}^{^{H}}(n,d)$ has all the properties stated in the corollary.\epr

\brem It is however not clear whether the moduli space of torsion-free Hitchin pairs ${\eM}_{_{S}}^{^{H}}(n,d)$ is even flat over $S$.  Of course, one does not expect good singularities on it. \erem

\section{The Gieseker-Hitchin map over a base} Let ${\ca}_{_S} \to S$ be the affine $S$--scheme representing the functor $T \mapsto \bigoplus_{i = 1}^{n} H^0(X_{_T}, \cl_{{_{X_{_T}}}}^i)$, where $X_{_T} = X \times_{_S} T$. Recall that we have chosen $\cl$ to be {\em sufficiently very ample} on $X/S$ so that the higher direct images $R^1f_{_*}(\cl^i)$ are zero. Then, it is not hard to see by the projection formula this is an affine \text{$S$-scheme} representing the functor $T \mapsto \bigoplus_{i = 1}^{n}  H^0(S,f_{_*}(\cl_{{_{X_{_T}}}}^i))$. This is the relative version of the space of characteristic polynomials. The points of this space ${\ca}_{_S}$ are polynomials $t^n + \sum_{i=1}^{n} q_i.t^{n-i}$, with $q_i \in H^0(S,f_{_*}(\cl^{^i}))$. 

Equivalently, a point of ${\ca}_{_S}$ can be viewed as a polynomial written $t^n + \sum_{i=1}^{n} q_i.t^{n-i}$, with $q_i \in H^0(X, \cl^{^i})$.

We observe that for each modification $X_{_T}^{^{(\text{mod})}}$ we have a canonical identification:
\beqa
{\ca}_{_S}(T) = \bigoplus_{i = 1}^{n} H^0(X_{_T}^{^{(\text{mod})}}, \cl_{_{\text{mod}}}^{^i})
\eeqa
\brem\label{undefinedness of Hitchin map} {\it Obstructions to defining the Hitchin map.} \begin{itemize}
\item Let $(E,\theta)$ be a torsion-free Hitchin pair on the surface $X$ and let $\eE$ be the pure sheaf on $Z$ which corresponds to $(E,\theta)$. Because the generic fibre is smooth and projective curve one can define the Hitchin map {\it $\sf h$} as is shown in \cite{sim2}. 
\item However, in general the Hitchin map  $\sf h$ does not extend in a well-defined manner to  the moduli functor of torsion-free Hitchin pairs on $X$; for instance,  if $T$ is an arbitrary \text{$S$-scheme}, then it is not clear why the values of the characteristic polynomial coincide when we approach the node on the closed fibre through distinct curves on $T$. 
\item Secondly, for an arbitrary \text{$S$-scheme} $T$, the singularities of $X_{_T}$ are no longer of the $A$-type. 
\item One of the important points of this paper is that the Hitchin map is well-defined on the Gieseker-Hitchin functor (Definition \ref{hitchin map}); indeed, the Gieseker-Hitchin space {\it resolves} the rational morphism $\sf h$.
\end{itemize}  
\erem

\bdefe\label{hitchin map} Define the morphism:
\beqa
{\sf g}_{_S}:{\underline{\eG}}_{_{S}}^{^{H}}(n,d) \to {\ca}_{_S}
\eeqa
for each \text{$S$-scheme} $T$, ${\sf g}_{_T}\big[X_{_T}^{^{(\text{mod})}},(V_{_T},\phi_{_T})\big] =  (q_1(\phi_{_T}), \ldots, q_{n-1}(\phi_{_T}))$.
\edefe 
\brem The the morphism ${\sf g}_{_S}$ respects the {\em equivalence} (Definition \ref{equivalence}) of families. The action of the automorphism $g$ (which leaves the end points $p_{_1}$ and $p_{_2}$ fixed does not affect the gluing, while on the bundle $V|_{_R}$ restricted to the tree $R$, $g$ acts as an automorphism and the Higgs structure $\phi$ is simply an endomorphism ($\cl_{_{\text{mod}}}$ is trivial on $R$). Thus, the Higgs structure is simply conjugated by the automorphism $g$ and hence the characteristic polynomial is well-defined. \erem

\subsubsection{The Spectral variety}

\bprop\label{weakproper} Let $X \to S$ be as before with a singular fibre $C$. Let $T \to S$ be a smooth curve over $S$ equipped with a marked point $\tau \in T$ over the closed point $s \in S$. Let $U = T - \{\tau\}$ and suppose that $(E_{_u},\theta_{_u})$ is a family of stable Hitchin pairs on the family of smooth curves $X_{_u}$ parametrized by $U$. Assume that the characteristic polynomial ${\sf h}(\theta_{_u})$ has a limit ${\sf h}_o$ in ${\ca}_{_T}$. Then, $\lim_{u \to \tau}(E_{_u},\theta_{_u})$ exists as a stable torsion-free Hitchin pair with ${\sf h}(\theta_{_\tau}) = {\sf h}_o$. \eprop
\bpr The proof follows \cite{sim2} closely, but we need to keep in mind two points; firstly, the closed fibre is singular and secondly the Hitchin map is not defined on the moduli space of torsion-free Hitchin pairs. 

Let $\eta_i(u) = q_i({\sf h}(\theta_{_u}))$, the coefficients of the characteristic polynomial. By assumption, these functions extends to the whole of $T$. Now consider the function $c:T \times_{_S} \cl \to \cl^{n}$ defined by $c(z,t) = t^n + \sum_{i=1}^{n} \eta_i(z).t^{n-i}$. By Cayley-Hamilton theorem, one knows that $c(u,\theta_{_u}) = 0$ for all $u \in U$.

Observe that the family $(E_{_u},\theta_{_u})$ gives a morphism $g':U \to {\eM}_{_S}^{^{H}}$. Using the compactification ${\eM}_{_S}^{^{H}} \subset M(\co_{_Z}, P_k)$ (see Lemma \ref{keycorresp}), we extend the map $g'$ to a map $g:T \to M(\co_{_Z}, P_k)$. By possibly going to a finite covering of $U$ and using the GIT construction of $M(\co_{_Z}, P_k)$, via $g$, we get a family $\eE$ on $Z \times_S T$, such for each $u \in U$, this family $\eE_{_u}$ gives a point of  ${\eM}_{_S}^{^{H}}(T)$. 

By Lemma \ref{keycorresp}, we have $Supp(\eE_{_u}) \cap D = \emptyset$ for each $u \in U$, i.e $Supp(\eE|_{_U}) \subset U \times_{_S} \cl$. Since $c(u,\theta_{_u}) = 0$ for all $u \in U$ the support is contained in the zero scheme $Z(c) \subset T \times_{_S} \cl$. 

Let $\eE_{_\tau} = \lim_{u \to \tau} \eE_{_u}$. By flatness, the support of $\eE_{_\tau}$ is contained in the closure of $Supp(\eE|_{_U})$ in $T \times_{_S} Z$. But this is contained in the closed subscheme $Z(c)$ and hence contained in $T \times_{_S} \cl$. In other words, the limiting sheaf 
$\eE_{_\tau}$ also has the property that $Supp(\eE_{_\tau}) \cap D = \emptyset$. Thus, by Lemma \ref{keycorresp}, we get a limiting Hitchin pair $(E_{_\tau}, \theta_{_\tau})$.\epr
\brem\label{on support} From the proof of Proposition \ref{weakproper} it follows that, if for some $T$, and a family $(E_{_T}, \theta_{_T})$ of torsion-free Hitchin pairs, if the characteristic polynomial of $\theta_{_T}$ is defined as a $T$-valued point of $\ca_{_S}$, then, one can consider the notion of a spectral scheme associated to a family of Hitchin pairs $(E_{_T}, \theta_{_T})$.
Once the characteristic polynomial, say $u$, is fixed in $\ca_{_S}(T)$, it defines a closed subscheme  $Y_{_u} \subset \cl$ and the family ${\eE}_{_T}$ of pure sheaves on $Z_{_T}$ which corresponds to $(E_{_T}, \theta_{_T})$ is {\em set-theoretically} supported on $Y_{_u}$. If moreover, $Y_{_u}$ is {\em reduced and irreducible}, it is precisely the {\em scheme theoretic support} of the pure sheaf ${\eE}_{_T}$ on $Z_{_T}$.\erem 

\bth\label{Hitchinproper} The Hitchin map \eqref{hitchin map}, ${\sf g}_{_S}:{\eG}_{_{S}}^{^{H}}(n,d) \to {\ca}_{_S}$ is proper over $S$. \eeth
\bpr Observe that over the generic point $\zeta \in S$, the conditions of  Lemma \ref{quasiprojmorphs} hold good by the classical properness of the Hitchin map (\cite[Theorem 6.11]{sim2}). Thus, we need to check only the {\it horizontal properness} of the Hitchin map. 

Let $T \to S$ be a smooth curve over $S$ equipped with a marked point $\tau \in T$ over $s \in S$. Let $U = T - \{\tau\}$ which maps to the generic point of $S$. 

Suppose that we are given a morphism $\alpha:U \to {\eG}_{_{S}}^{^{H}}(n,d)$ such that ${\sf g}_{_S} \circ \alpha$ extends to the whole of $T$. We need to check that $\alpha$ itself extends to $T$. Now since over the subset $U$ the surface is isomorphic to $X_{_U}$, the map $\alpha$ gives a family of stable Hitchin pairs $(E_{_u}, \phi_{_u})$ on $X_{_U}$. Observe that since the line bundle $\cl_{_{\text{mod}}}$ is a pull-back of $\cl$ from $X$, the characteristic polynomial $\lim_{u \to \tau} {\sf g}_{_U}(E_{_u}, \phi_{_u})$ exists in ${\ca}_{_S}(T)$ (since ${\sf g}_{_S} \circ \alpha$ extends to the whole of $T$). 

By Proposition \ref{weakproper}, the family extends to a torsion-free stable Hitchin pair $(E_{_\tau}, \phi_{_\tau})$. Now applying Corollary \ref{giestotfproper1}
, we see that $\alpha$ extends to $T$ completing the proof.
\epr
In summary we have the following theorem:
\bth\label{summary} Let $X/S$ be as before.
\begin{enumerate} 
\item We have a quasi-projective \text{$S$-scheme} ${\eG}_{_{S}}^{^{H}}(n,d)$ which is flat over $S$ which is regular as a scheme over $k$. 
\item The generic fibre $({\eG}_{_{S}}^{^{H}}(n,d))_{_\zeta}$ is the classical Hitchin moduli space of stable Hitchin pairs on $X_{_\zeta}$ of rank $n$ and degree $d$, and the closed fibre $({\eG}_{_{S}}^{^{H}}(n,d))_{_s}$ is a divisor with analytic normal crossing singularities.
\item We have a natural Hitchin morphism ${\sf g}_{_S}:{\eG}_{_{S}}^{^{H}}(n,d) \to {\ca}_{_S}$ which is proper over $S$ and extends the classical Hitchin morphism on $({\eG}_{_{S}}^{^{H}}(n,d))_{_\zeta}$.
\end{enumerate} \eeth

Recall the $S$-morphism $\nu_{_*}:{\eG}_{_{S}}^{^{H}}(n,d) \to {\eM}_{_{S}}^{^{H}}(n,d)$, which is an isomorphism over the generic point $\zeta \in S$.

\bprop\label{descentofhitchinmap} The Hitchin map ${\sf g}_{_S}:{\eG}_{_{S}}^{^{H}}(n,d) \to {\ca}_{_S}$ descends to a well-defined {\em set-theoretic} map  on the image $Im(\nu_{_*})$
and we have a diagram: 
\beqa\label{hitchintriangle2}
\xymatrix{
{\eG}_{_{S}}^{^{H}}(n,d) \ar[d]_{{\sf g}_{_S}} \ar[r]^{\nu_{_*}} & Im(\nu_{_*}) \ar[ld]^{{\sf h}_{_S}} \\
{\ca}_{_S}  
}
\eeqa 
\eprop
\bpr As $\nu_{_*}$ is an isomorphism outside the closed point on $S$, it is enough to check the statement on the closed fibre $X_{_s} = C$. Let $(E,\theta)$ is a stable Hitchin pair on $C$ so that $\theta:E \to E \otimes \cl$. Let the local type of $E$ at the maximal ideal be $\co^a \oplus {\mathfrak m}^b$; choose any $(V,\phi)$,  on the semistable curve $C^{(b)}$ such that $\nu_{_*}(V,\phi) = (E,\theta)$ (by assumption $(E,\theta) \in Im(\nu_{_*})$). Note that $\phi:V \to V \otimes \nu^*(\cl)$. 

As $\nu^*(\cl)|_{_R} = \co_{_R}$, it follows that the coefficients of the characteristic polynomial $q_i(\phi)$ give sections of $\cl_{_{\tilde{C}}}^i$ with an identification at the points $p_1, p_2$. This gives sections of $\cl^i$ on $C$. Clearly these are independent of the $(V,\phi)$ chosen on a semistable curve above $C$ since the characteristic polynomial is defined in terms of the restrictions of $V$ to $\tilde{C}$, and by \cite[Remark 4(i), page 176]{ns2}, the direct image $\nu_{_*}(V) = E$ determines $V|_{_{\tilde{C}}}$.  Hence they define the  characteristic polynomial of $(E,\theta)$.\epr
\section{The Hitchin fibre} 

 Let $X \to S$ be as before a fibered surface with a singular fibre $C$ which is irreducible with a single node. Also, $X$ is regular over $k$. By Theorem \ref{Hitchinproper}  the Hitchin morphism is well-defined and proper on the Gieseker-Hitchin scheme ${\eG}_{_{S}}^{^{H}}(n,d)$ and it is {\em not} well-defined on the space of Hitchin pairs ${\eM}_{_{S}}^{^{H}}(n,d)$.  We recall briefly well-known facts from Hitchin's work and highlight facts which are specially relevant for our discussion. We make the following general observations first. 

\blem\label{generalspectralsurface} Let $p:X \to S$ be as before. Let $\cl$ be a line bundle such that $\cl$ is {\em relatively very ample}.  For a general section $\xi:S \to  \oplus_{i=1}^{n} p_{_*}(\cl^i)$, we get a spectral surface $\psi_{_\xi}: Y_{_\xi} \to X$,  which is a ramified $n$-sheeted cover over $X$ such that it is unramified over the nodes of the special fibre $X_{_s}$.\elem 
\bpr The proof is essentially from  \cite[page 172]{bnr}. We quickly recall that the construction of the spectral surface $Y_{_\xi}$ with the associated properties. Let $W = {\underline {\spec}}(Sym(\cl^*))$ be the total space of the line bundle $\cl$ and $\varrho:W \to X$ the projection. The pull-back $\varrho^*(\cl)$ then gets a tautological section $x$. Now take the sections of $\varrho^*(\cl^n)$ on $W$ of the form 
\beqa
\xi = x^n + \sum_{i=1}^{n} \varrho^*(\xi_{i}).x^{n-i}
\eeqa
for $\xi_i:S \to  p_{_*}(\cl^i)$. The zero scheme $Y_{_\xi}$ of the sections $\xi$ of $\varrho^*(\cl^n)$ will be the spectral surface we desire. 

We need to show the existence of one such which satisfies the properties of the Proposition. Consider the embedding $X \hra {\bp}^{^N} \times S$ given by the line bundle $\cl$. Restrict this to the closed fibre $X_{_s} \subset {\bp}^{^N} \times \{s\}$ and choose a section $h \in H^{^0}(X_{_s}, \cl_{_s}^{^{n}})$ such that the hyperplane in $\bp(H^{^0}(X_{_s}, \cl_{_s}^n))$ defined by $h$ has the following properties:
\begin{enumerate}
\item $h$ does not meet the nodes on $X_{_s}$. Because $\cl_{_s}$ is very ample, this is always possible. In other words, 
\beqa
(h)_{_{0}} \cap \{nodes~~of~~ X_{_s}\} = \emptyset
\eeqa
\item $h$ does not have multiple zeros
\item $Y_{_h}$ is irreducible.
\end{enumerate}
Again, we can choose a section $\xi_{_n}:S \to p_{_*}(\cl^n)$ such that $\xi_{_n}(s) = h$ and by choice, we see that:
\beqa
(\xi_{_n}(s))_{_0} \cap \{nodes~~of~~ X_{_s}\} = \emptyset
\eeqa

We now consider the special section $\xi = x^n + \varrho^*(\xi_{_n})$. By the properties of the section $h$, it follows that the spectral curve $Y_{_h}$ defined by $x^n + \varrho^*(h)$ is {\em smooth except for nodal singularities}, with exactly $n$-nodes over each node of $X_{_s}$. The smoothness is the consequence of the fact that $h$ has no multiple zeroes and the discriminant of this polynomial is (upto a sign) simply $n^n.\varrho^*(h)^{n-1} \in H^0(X_{_s},\cl^{n(n-1)})$. 

By the choice of $\xi_{_n}$ as a generization of $h$ and by the {\it openness} of smoothness, the generic spectral curve  $Y_{_{\xi_{_{_n}}(\zeta)}} \to X_{_\zeta}$ defined by the section $\xi_{_n}(\zeta)$ is smooth. Thus the spectral surface $\psi_{_\xi}:Y_{_\xi} \to X$ has all the properties that was claimed. This shows that the set of spectral surfaces with these properties is non-empty. It is clearly open and the result follows.\epr

We now define the  {\em genericity} condition to analyze the general Hitchin fibre, which by Lemma \ref{generalspectralsurface} is non-empty.
\bdefe\label{genericity} 
We define the subspace:
\beqa\label{urs}
\Gamma({\cA}_{_S}^{^{ur}}):= \{sections~~\xi:S \to {\ca}_{_S}\mid \xi~~ is~as~in~Lemma~\ref{generalspectralsurface} \}
\eeqa
to be the set of {\em general} sections of ${\ca}_{_S} \to S$. \edefe
\brem  Let $\xi:S \to {\ca}_{_S}$ be a general section. Then there is a canonically defined spectral fibered surface $Y_{_\xi}$ over $S$ together with a covering $S$-morphism $\psi_{_\xi}:Y_{_\xi} \to X$; over the generic point $\zeta \in S$, $(\psi_{_\xi})_{_K}:(Y_{_\xi})_{_K} \to X_{_K}$ it is the classical spectral curve which is {\em smooth and irreducible} and over the closed fibre $s \in S$ when for instance the fibre $X_{_s} = C$ is a nodal curve with a single node  the spectral cover $\psi_{_u}:Y_{_{u}} \to C$, for $u = \xi(s)$ 
is as in Figure 3 where the fibre $Y_{_{u}}$ is a {\em vine curve}. \erem 
\begin{figure}
 \begin{tikzpicture}[scale=0.3]
\draw plot[smooth cycle] coordinates {(-5,1)     (-3.5,0)  (-1, 2)  (0.75,1) (-1,0) (-3.5, 2) };
\draw plot[smooth] coordinates {(1.5, 2) (0.75,1) (1.5, 0)};
 \draw   [<-, thick] (-8,1) -- (-6.55,1);
  \draw   (-11.5,0) .. controls (-3.5,4) and (-17,4) .. (-8.5,0);
\end{tikzpicture}
\caption{{The spectral cover $C \leftarrow Y_{_{u}}$}}
\label{fig3}
\end{figure}
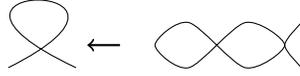

\brem\label{spectralproperty} Let $\xi$ be as above and $u = \xi(s)$. Let $\psi_{_u}:Y_{_u} \to C$ be the spectral curve over the closed point $s \in S$. The direct image $(\psi_{_u})_*(\co_{_{Y_{_u}}}) = \co_{_{Y}} \oplus \cl^* \ldots \oplus (\cl^*)^{^n}$ and the genus $g_{_{Y_{_u}}}$ of ${Y_{_u}}$ is equal to $n(g - 1) + 1 + deg(\cl).\frac{n(n-1)}{2}$, where $g$ is the genus of $C$.

By  \cite{hitchin} (or \cite[Proposition 3.6]{bnr}), we have a bijective correspondence between torsion-free sheaves $\eta$ on $Y_{_\xi}$ of rank $1,$ relative degree $\delta$ and families of {\em stable} torsion-free Hitchin pairs $(E,\theta)$ on $X$, where $E$ is torsion-free of rank $n$ and degree $d = \delta - deg(\cl)\frac{n(n-1)}{2}$ and $\theta:E \to E \otimes \cl$ is a homomorphism with characteristic coefficients $\xi_i:S \to p_{_*}(\cl^i)$. The correspondence is given as follows. Let $W = {\underline {\spec}}(Sym(\cl^*))$ be the total space of the line bundle $\cl$. Recall the diagram:
\beqa\label{spectraltriangle}
\xymatrix{
Y_{_\xi} \ar[d]_{\psi_{_\xi}} \ar[r]^{\subset}& W \ar[ld]^{\pi} \\
X  
}
\eeqa
The line bundle $\pi^*(\cl)$ has a tautological section $t$ which induces the canonical map
\beqa\label{spectralhitchinmap}
\eta \stackrel{1 \otimes t}\longrightarrow \eta \otimes \psi_{_\xi}^*(\cl)
\eeqa
Pushing this down gives the map $\theta:E \to E \otimes \cl$, where $E:= \pi_*(\eta)$. The correspondence sends $\eta \mapsto (E, \theta)$. The stability of the Hitchin pair $(E,\theta)$ is easily checked.
\erem

If $\xi:S \to {\ca}_{_S}$ is any section then we get a spectral surface $Y_{_\xi} \subset W = {\underline {\spec}}(Sym(\cl^*)) \subset Z$, where $Z = {\bp}(\cl^* \oplus \co_{_X})$, is the projective completion of the total space $W$ of $\cl$ as a scheme over $S$ (see Lemma \ref{keycorresp}).
\bprop\label{picaspure} Let $\xi:S \to {\ca}_{_S}$ be a general section. Then the compactified Picard variety ${P_{_{\delta,Y_{_\xi}}}}$  of spectral fibered surface $Y_{_\xi} \subset Z$ can be canonically identified with the subscheme of the moduli space ${M}(\co_{_Z}, P_k)$ of pure sheaves $\eE$ on $Z$ such that the scheme theoretic support $Supp(\eE) = Y_{_\xi}$. \eprop
\bpr Recall (Lemma \ref{keycorresp}) that a family  $(E,\theta)$ of stable torsion-free Hitchin pairs  on $X$, canonically defines a family of stable pure sheaves $\eE$ on the scheme $Z$ over $S$. As $\xi$ is generic by Remark \ref{spectralproperty}, rank $1$ torsion-free sheaves of relative degree $\delta$ on $Y_{_\xi}$ give points of ${\eM}_{_{S}}^{^{H}}(n,d)$ or equivalently pure sheaves on $Z$. Again by the genericity of $\xi$,  the scheme $Y_{_\xi}$ is reduced and irreducible and hence by Remark \ref{on support}, these pure sheaves $\eE$ have scheme theoretic support $Supp(\eE) = Y_{_\xi}$. i.e the compactified Picard variety ${P_{_{\delta,Y_{_\xi}}}}$  (cf. Caporaso \cite{caporaso}) gets realized as a subscheme of ${M}(\co_{_Z}, P_k)$ which parametrizes $\tt p$--semistable pure sheaves on $Z$ with fixed Hilbert polynomial. \epr

Recall that the scheme structure on ${\eM}_{_{S}}^{^{H}}(n,d)$ was realized as an open subscheme of ${M}(\co_{_Z}, P_k)$ which parametrizes $\tt p$--semistable pure sheaves on $Z$ with fixed Hilbert polynomial.

Thus by Proposition \ref{picaspure}, we have an inclusion:
\beqa\label{obs}
{P_{_{\delta,Y_{_\xi}}}} \subset {\eM}_{_{S}}^{^{H}}(n,d)
\eeqa

\bth Let $\xi$ be a general section as above and let the fibre of ${\sf g}_{_S}$ over $\xi$ be denoted by ${\sf g}_{_S}^{-1}(\xi)$. Then we have a proper birational morphism of \text{$S$-scheme}s:
\beqa
{\sf g}_{_S}^{-1}(\xi) \to {P_{_{\delta,Y_{_\xi}}}}
\eeqa
which is an isomorphism over $\spec~K$; more precisely, it coincides over $\spec~K$ with the classical Hitchin isomorphism of the Hitchin fibre with the Jacobian of a smooth spectral curve $(Y_{_{\xi}})_{_K}$.
\eeth 
\bpr By Proposition \ref{descentofhitchinmap} we see that for every \text{$S$-scheme} $T$,\small
\beqa
\nu_{_*}({\sf g}_{_S}^{-1}(\xi))(T) = \{(E,\theta)\in {\eM}_{_{S}}^{^{H}}(n,d)(T) \mid Supp(\eE) = Y_{_\xi} \times_{_S} T \}
\eeqa\normalsize
and by the observation \eqref{obs}, we get a proper birational surjective morphism:
\beqa
\nu_{_*}:{\sf g}_{_S}^{-1}(\xi) \to {P_{_{\delta,Y_{_\xi}}}}
\eeqa
\epr
\brem By Zariski's Main theorem, since ${P_{_{\delta,Y_{_\xi}}}}$ is normal (\cite{caporaso}), $\nu_{_*}$ has connected fibres. Since the morphism $\nu_{_*}$ is an isomorphism over the generic fibre, we need to look closely on the phenomenon over the closed fibre i.e the nodal curve $C$. \erem

\section{Geometry of the degenerate Hitchin fibre}

The aim of this section is to give a description of the geometry of the Hitchin fibre and prove a statement which can be described as a {\it quasi-abelianization}
of the moduli space of Hitchin pairs.
\subsubsection{A review of the compactified Picard variety}
We begin by a variation in the description of the compactification of the Picard variety of a stable curve. For the sake of simplicity we work with an irreducible {\em vine curve} $Y$ with $n$-nodes (which occurs as our spectral curve) and take a re-look at the compactification of the Picard variety of $Y$. 

Recall that since $Y$ is irreducible, there is a natural choice of the compactification. In \cite{alex}, we find a comparison of various approaches to the compactification, beginning with the one by Oda-Seshadri \cite{os}, Caporaso \cite{caporaso} and Simpson \cite{sim1}. In fact, all three approaches give the same object.  Recall that in \cite{os}, the compactification is described as a moduli of torsion-free sheaves on the curve $Y$ with fixed slope while in \cite{caporaso}, following Gieseker, the description is in terms of embeddings of semistable curves stably equivalent to the curve $Y$. 

The description we wish to give here is closer in spirit to the one in \cite{caporaso} and comes from the paper of Nagaraj-Seshadri \cite{ns2}. Following the approach in \cite{ns2} (see also \cite[page 15]{ictp}) we realize the compactification of $Pic~Y$,  as Gieseker line bundles on a {\it ladder} curve semistably equivalent to $Y$ (see Proposition \ref{oscaporaso}). This approach is essential in our description of the geometry of the Gieseker-Hitchin fibre which reveals more interesting phenomena and a new compactification of the Picard variety (see Theorem \ref{relative abelianization} and Remark \ref{morecomparison} below for details).

\begin{figure}[htbp]
\begin{center}
\includegraphics[scale=0.80]{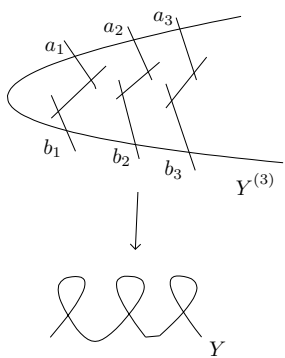} 
\caption{{\bf The contraction of $R^{^{(3)}}$'s}}
\label{Gieseker3}
\end{center}
\end{figure}

Let $Y^{^{(\ell)}}$ be the semistable curve obtained from $Y$ by attaching trees of length $\ell$ to the normalization of $Y$ at the points $a_i,b_i, i = 1, \ldots, n$ over the $n$-nodes on $Y$. Let $p:Y^{^{(\ell)}} \to Y$ be the morphism contracting the trees $R^{^{(\ell)}}(a_i,b_i)$ joining the pairs of points $a_i,b_i, i = 1,2, \ldots, n$ (see Figure \ref{Gieseker3}).

\bdefe\label{Giesekerlinebundle} A quasi-Gieseker line bundle $\eN$ on $Y^{^{(\ell)}}$  is  a line bundle such that on each tree $\{R^{^{(\ell)}}(a_i,b_i)\}_{i=1}^n,$ we have the line bundle $\eN|_{_{{R^{^{(\ell)}}}(a_i,b_i)}}$ is {\em standard} and 
such that the direct image $p_{_*}(\eN)$ is torsion free on $Y$.
\edefe
\brem The key difference between this definition of a quasi-Gieseker line bundle and a Gieseker vector bundle earlier in Definition \ref{Giesekervectorbundle} is that in Definition \ref{Giesekervectorbundle}  we impose the {\em strict standardness}. In the quasi-Gieseker line bundle we allow the line bundle to be {\em trivial} on some $\bp^{^1}$'s in some trees while in the Gieseker vector bundle definition this is not allowed. Also, the considerations in Definition \ref{Giesekervectorbundle} were for a nodal curve with a single node.\erem

\bdefe\label{auto1} A {\em vertical automorphism} $\sigma:Y^{^{(\ell)}} \to Y^{^{(\ell)}}$ which leaves the end-points $a_i$ and $b_i$ fixed is an automorphism:
\beqa\label{auto2}
\xymatrix{
{Y^{^{(\ell)}}} \ar[dr]_{p} \ar[rr]^{\sigma}& & Y^{^{(\ell)}}\ar[dl]^{p} \\
& Y  &
}
\eeqa
which acts on the trees $R^{^{(\ell)}}(a_i,b_i)$ leaving the end-points $a_i$ and $b_i$ fixed.  The group of vertical automorphisms leaving the end-points of the trees fixed will be denoted by $Aut_{_{epf}}(Y^{^{(\ell)}})$. \edefe

\brem\label{autoeg} Take the case when $\ell = 1$. As an automorphism $\lambda$ of a tree $R^{^{(1)}}$  which fixes two points is the multiplicative group $\bg_{_m}$,  we see that $\sigma = \big(\lambda_{_1}, \lambda_{_2}, \ldots, \lambda_{_n}\big) \in  \bg_{_m} \times \ldots \times \bg_{_m}$,
i.e 
\beqa
Aut_{_{epf}}(Y^{^{(1)}}) \simeq \overbrace{\bg_{_m} \times \ldots \times \bg_{_m}}^{n}.
\eeqa 
Each component $\lambda_i$ act trivially on $\co_{_{R^{^{(1)}}}}$, while on $\co_{_{R^{^{(1)}}}}(1)$ it acts as multiplication by $\lambda_i$ on the fibre $\co_{_{R^{^{(1)}}}}(1)_{_{a_i}}$ and by $\lambda_i^{^{-1}}$ on $\co_{_{R^{^{(1)}}}}(1)_{_{b_i}}$ for each $i$. In other words, we have a canonical action of such automorphisms on the set of quasi-Gieseker line bundles on the curve $Y^{^{(1)}}$.

\begin{figure}[htbp]
\begin{center}
\includegraphics[scale=0.80]{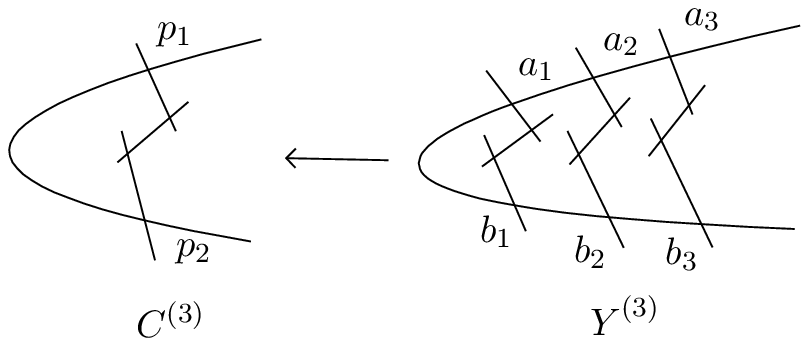} 
\caption{{\bf The covering $Y^{^{(3)}} \to C^{^{(3)}}$}}
\label{Gieseker4}
\end{center}
\end{figure}

Similarly we see that:
\beqa
Aut_{_{epf}}(Y^{^{(\ell)}}) \simeq \overbrace{\bg_{_m}^{^\ell} \times \ldots \times \bg_{_m}^{^\ell}}^{n}
\eeqa
In the earlier setting of Definition \ref{equivalence}, when we work with the curve $C^{^{(\ell)}}$, with a single tree $R^{^{(\ell)}}(p_{_1},p_{_2})$ joining $p_{_1}$ and $p_{_2}$, the group of vertical automorphisms is simply $Aut_{_{epf}}(C^{^{(\ell)}}) \simeq \bg_{_m}^{^\ell}$. The equivalence of Gieseker-Hitchin pairs is defined via the orbits of this group (see Figure \ref{Gieseker4}).
\erem
\bdefe Let \small
\beqa
{\eG}_{_{Y}}(1,\delta):= \frac{\{quasi-Gieseker~~line ~~bundles~~on~~Y^{^{(1)}}~of~deg~\delta\}}{Aut_{_{epf}}(Y^{^{(1)}})}
\eeqa\normalsize
the isomorphism classes of quasi-Gieseker line bundles modulo  the action of $Aut_{_{epf}}(Y^{^{(1)}})$ as in Definition \ref{auto1}.\edefe

Following the strategy of \cite{ns2} or Section 6 of this paper, one can give a natural scheme structure to the set ${\eG}_{_{Y}}(1,\delta)$ which we will call the moduli space of quasi-Gieseker line bundles on $Y^{^{(1)}}$.
\bprop\label{oscaporaso} The compactified Picard variety ${P_{_{\delta,Y}}}$ of $Y$ is isomorphic to the moduli space ${\eG}_{_{Y}}(1,\delta)$ of quasi-Gieseker line bundles on $Y^{^{(1)}}$. The isomorphism is induced by the direct image morphism $p_{_*}$. \eprop  
\bpr This result is proven by Pandharipande \cite{pandharipande}. In the context of the present paper, the proof can be given as in \cite[Theorem 2, page 196]{ns2} (see also Section 6 above), where this isomorphism is shown more generally for the case when the rank and degree are coprime, except that in \cite{ns2}, the case is when the curve $Y$ has a single node. The generalization to all stable curves has been carried out in \cite{schmitt}.\epr
\brem When the number of nodes is strictly bigger than $1$, the singularities of the compactified Picard variety is a {\em product of normal crossing singularities} and therefore {\em not} a normal crossing singularities (cf. \cite[Page 595]{caporaso}, \cite[Page 262, I]{ictp}). This is erroneously written in \cite[Theorem 3.3.1]{schmitt}. \erem

\subsubsection{Towards the geometry of the Gieseker-Hitchin fibre}

\begin{figure}[htbp]
\begin{center}
\includegraphics[scale=0.80]{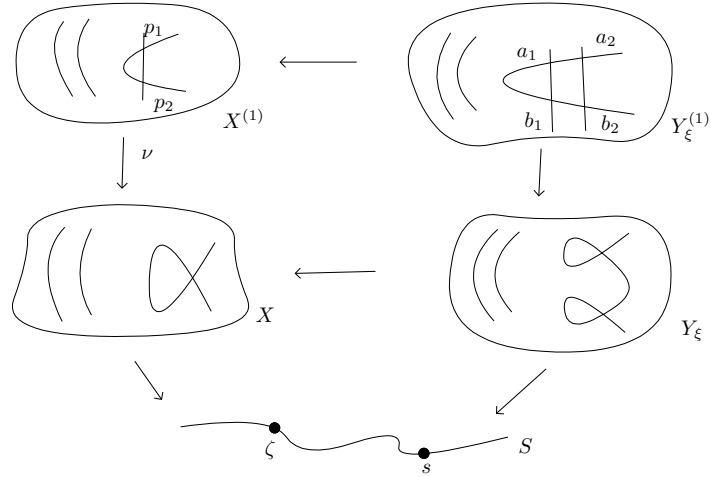} 
\caption{{\bf The spectral picture}}
\label{The spectral picture}
\end{center}
\end{figure}

\blem Let $\xi \in \Gamma({\cA}_{_S}^{^{ur}})$  as in Definition \ref{genericity}. Then, for each $\ell$, the section $\xi$ defines a ({\it spectral}) covering surface $\psi_{_\xi}^{(\ell)}:Y^{(\ell)}_{_\xi} \to X^{(\ell)}$, where $Y^{(\ell)}_{_\xi}$ is a (semistable) modification of $Y_{_\xi}$. Over the closed point $s \in S$ the fibre $Y^{(\ell)}_{_{\xi,s}}$ is  a curve with $n$ pairs of marked points each pair being joined by a  tree $R^{^{(\ell)}}$ and the covering morphism $Y^{(\ell)}_{_{\xi,s}} \to C^{(\ell)}$ as in Figure \ref{Gieseker4}.
\elem
\bpr By the choice of $\xi$, we have a covering morphism $\psi_{_\xi}:Y_{_\xi} \to X$ which has the good properties given by Lemma \ref{generalspectralsurface}. Now by \cite[Corollary 7.15, Chapter II]{hartshorne}, we have a natural morphism  
$\psi_{_\xi}^{(\ell)}:Y^{(\ell)}_{_\xi} \to X^{(\ell)}$ and a diagram:
\beqa\label{keydiag1}
\xymatrix{
X^{(\ell)}  \ar[d]_{\nu} &
\ar[l]_{{\psi_{_\xi}^{(\ell)}}} Y^{(\ell)}_{_\xi} \ar[d]^{p} \\
X  & \ar[l]^{\psi_{_\xi}} Y_{_\xi} 
}
\eeqa
where $p$ and $\nu$ contracts the $R^{^{(\ell)}}$'s to the respective stable curves (see Figure \ref{The spectral picture}). 

Consider the morphism $\nu:X^{(\ell)} \to X$ and $\cl_{_\ell}$ the pull-back $\nu^{^*}(\cl)$. The generic section $\xi$ as in Lemma \ref{generalspectralsurface} pulls back to give a section $\xi_{_\ell}:S \to  \oplus_{i=1}^{n} p_{_*}(\cl_{_\ell}^i)$. 
Let $W_{_\ell} = {\underline{\spec}}(Sym(\cl_{_\ell}^*))$ be the total space of the line bundle $\cl_{_\ell}$. As in Lemma \ref{generalspectralsurface}, we can take the spectral surface $Y_{_{\xi_{_\ell}}}$ defined by the section $\xi_{_\ell}$ as a subscheme of $W_{_\ell}$. We have a canonical diagram:
\beqa\label{keydiag1.1}
\xymatrix{
X^{(\ell)}  \ar[d]_{\nu} &
\ar[l]_{{\psi_{_{\xi_{_\ell}}}}} Y_{_{\xi_{_\ell}}} \ar[d]^{p} \\
X  & \ar[l]^{\psi_{_\xi}} Y_{_\xi} 
}
\eeqa
By the universal property of blow-ups, it is easy to see that $Y_{_{\xi_{_\ell}}} \simeq Y^{(\ell)}_{_\xi}$. The remaining claims in the Lemma are easily established. 
\epr

\brem\label{automorphisms} Let $T$ be a \text{$S$-scheme} and let $\sigma:Y^{(\ell)}_{_{\xi,T}} \to Y^{(\ell)}_{_{\xi,T}}$ be an automorphism as follows:
\beqa\label{auto}
\xymatrix{
{Y^{(\ell)}_{_{\xi,T}}} \ar[dr]_{p} \ar[rr]^{\sigma}& & Y^{(\ell)}_{_{\xi,T}}\ar[dl]^{p} \\
& Y_{_{\xi,T}}  &
}
\eeqa
In other words, $\sigma$ is essentially given by an automorphism over a point $t \in T$ above the closed point $s \in S$. i.e., the entire information is given by the following diagram: 
\beqa\label{autoovert}
\xymatrix{
{Y^{(\ell)}_{_{\xi,t}}} \ar[dr]_{p_{_t}} \ar[rr]^{\sigma_{_t}}& & Y^{(\ell)}_{_{\xi,t}}\ar[dl]^{p_{_t}} \\
& Y_{_{\xi,t}}  &
}
\eeqa
and by Figure \ref{The spectral picture} and Definition \ref{auto1} above, giving $\sigma_{_t} \in Aut_{_{epf}}(Y^{(\ell)}_{_{\xi,s}})$ is giving a tuple of $n$-automorphisms of the trees $R^{^{(\ell)}}(a_i,b_i)$'s joining the points $a_i$ and $b_i$ leaving the points $a_i$ and $b_i$ fixed. By Remark \ref{autoeg} we get a point of $\big(\lambda_{_1}, \lambda_{_2}, \ldots, \lambda_{_n}\big) \in  \bg_{_m}^{^\ell} \times \ldots \times \bg_{_m}^{^\ell}$.\erem

\brem \label{diagonallift} The morphism $Y^{(\ell)}_{_{\xi,t}} \to X_{_t}^{(\ell)}$ is a covering morphism which is unramified over the $R^{^{(\ell)}}(p_{_1}, p_{_2})$ (see Figure \ref{The spectral picture}). Hence the automorphism $\lambda$ of the $R^{^{(\ell)}}(p_{_1}, p_{_2})$  lifts to an automorphism of $Y^{(\ell)}_{_{\xi,t}}$ as a {\it diagonal} element $\lambda:= \big(\lambda, \lambda, \ldots, \lambda\big)$ in $\bg_{_m}^{^\ell} \times \ldots \times \bg_{_m}^{^\ell}$. 
\erem

\subsubsection{The stratification of the Gieseker-Hitchin spaces}
We work over the nodal curve $C$. Recall that the Gieseker-Hithin space ${\eG}_{_{C}}^{^{H}}(n,d)$ was constituted of Gieseker-Hitchin pairs $(V,\phi)$ on various curve $C^{^{(\ell)}}$ which were modifications of the nodal curve $C$. This can be expressed as giving the following stratification of the moduli spaces:
\beqa\label{ghstrata}
{\eG}_{_{C}}^{^{H}}(n,d) = \bigsqcup {\eG}_{_{C}}^{^{H}}(n,d)_{_{(\ell)}}
\eeqa
where 
\beqa\label{ghstrata1}
{\eG}_{_{C}}^{^{H}}(n,d)_{_{(\ell)}} := \frac{\{(V,\phi) \mid V~a~Gieseker~bundle~on~C^{^{(\ell)}}\}}{Aut_{_{epf}}(C^{^{(\ell)}})}
\eeqa
(see Remark \ref{autoeg}).

This stratification of the total space induces a stratification of the Gieseker-Hitchin fibre ${\sf g}_{_C}^{-1}(\xi_{_s})$ given as follows:
\beqa\label{ghfstrata}
{\sf g}_{_C}^{-1}(\xi_{_s}) = \bigsqcup {\sf g}_{_C}^{-1}(\xi_{_s})_{_{(\ell)}}
\eeqa
where
\beqa
{\sf g}_{_C}^{-1}(\xi_{_s})_{_{(\ell)}} := {\sf g}_{_C}^{-1}(\xi_{_s}) \cap {\eG}_{_{C}}^{^{H}}(n,d)_{_{(\ell)}}
\eeqa

Let $\eN$ be a  {\em a quasi-Gieseker line bundle} on $Y^{(\ell)}_{_{\xi,s}}$ in the sense of Definition \ref{Giesekerlinebundle}. For simplicity of notation, let us denote the morphism $\psi_{_{\xi,s}}^{(\ell)}:Y^{(\ell)}_{_{\xi,s}} \to C^{^{(\ell)}}$ for each $\ell$ simply by $\psi^{^{(\ell)}}$.
\bdefe Let $G^{^{(\ell)}}(1,\delta)$ denote the {\em set of isomorphism classes of quasi-Gieseker line bundles } with degree $\delta$ on the curve $Y^{(\ell)}_{_{\xi,s}}$. \edefe
\brem Note that we are not quotienting out by any equivalence of the automorphism action.\erem

\bprop\label{ss spectral curve} If $\eN \in G^{^{(\ell)}}(1,\delta)$ is a quasi-Gieseker line bundle on $Y^{(\ell)}_{_{\xi,s}}$, then, $\psi^{^{(\ell)}}_{_*}(\eN)$ is a stable Gieseker vector bundle on $C^{^{(\ell)}}$ equipped with a Higgs structure $\phi$ making $(\psi^{^{(\ell)}}_{_*}(\eN),\phi)$ a stable Gieseker-Hitchin pair.  Moreover, this induces an identification:
\beqa
{\sf g}_{_C}^{-1}(\xi_{_s})_{_{(\ell)}} \simeq G^{^{(\ell)}}(1,\delta)/\Delta_{_\ell}
\eeqa
where $\Delta_{_\ell}$ is the induced {\it diagonal} action of $\bg_{_m}^{^\ell}$ (see Remark \ref{diagonallift}).
\eprop
\bpr As in Remark \ref{spectralproperty}, we see that in the {\it spectral situation} such as $\psi_{_{\xi,s}}^{^{(\ell)}}:Y^{^{(\ell)}}_{_{\xi,s}} \to C^{^{(\ell)}}$ we have a corresponding diagram:  
\beqa\label{spectraltriangle1}
\xymatrix{
Y^{(\ell)}_{_{\xi,s}} \ar[d]_{\psi_{_\xi}^{(\ell)}} \ar[r]^{\subset}& W_{_\ell} \ar[ld]^{f_{_\ell}} \\
C^{^{(\ell)}}  
}
\eeqa
and it follows that if $\eN$ is any quasi-Gieseker line bundle on $Y^{^{(\ell)}}_{_{\xi,s}}$, then, $V = (\psi_{_\xi}^{(\ell)})_{_*}(\eN)$ is a vector bundle on $C^{^{(\ell)}}$ canonically equipped with a Higgs structure $\psi:V \to V \otimes \cl_{_\ell}$ since it arises from a spectral construction.

The commutativity of the diagram \eqref{keydiag1}, gives the identification
\beqa
(\psi_{_\xi})_{_*}(p_{_*}(\eN)) = \nu_{_*}(V)
\eeqa
from which we conclude that $\nu_{_*}(V)$ is the underlying torsion-free sheaf of a stable torsion-free Hitchin pair on $C$.  Hence by the definition of a stable Gieseker-Hitchin pair, it follows that $(V,\psi)$ obtained above is a {\em stable Gieseker-Hitchin pair}. 

The action of $Aut_{_{epf}}(C^{^{(\ell)}}) = \bg_{_m}^{^\ell}$ to determine the open startum of the Gieseker-Hitchin fibre (see \eqref{ghstrata1}) lifts to give the diagonal action on $G^{^{(\ell)}}(1,\delta)$, and we get the identification. \epr

The compactified Picard variety ${P_{_{\delta,Y_{_{\xi,s}}}}}$ of $Y_{_{\xi,s}}$ (which is an irreducible  vine curve has $n$-nodes) also has a stratification in terms of the complexity of the torsion-freeness of the sheaves. This can be given as follows:
\beqa\label{picstrata0}
{P_{_{\delta,Y_{_{\xi,s}}}}} = \bigsqcup {P_{_{\delta,Y_{_{\xi,s}}}}}(j)
\eeqa
where
\beqa\label{picstrata1}
{P_{_{\delta,Y_{_{\xi,s}}}}}(j):= \{\eta \mid \eta~is~\text{non-free}~at~exactly~j~nodes \}
\eeqa
In this description ${P_{_{\delta,Y_{_{\xi,s}}}}}(0)$ corresponds to the open subset of line bundles on $Y_{_{\xi,s}}$ of degree $\delta$.

\subsubsection{The big stratum}
Recall that we have a proper birational morphism $\nu_{_*}:{\sf g}_{_C}^{-1}(\xi_{_s}) \to {P_{_{\delta,Y_{_{\xi,s}}}}}$.
\bprop\label{opentoric} Let $\eta \in {P_{_{\delta,Y_{_{\xi,s}}}}}(n)$, in the {\it worst} stratum, i.e $\eta$ is given by the maximal ideal sheaf on each of the $n$-nodes.  The part of the fibre $\nu_{_*}^{^{-1}}(\eta)$ in the open stratum ${\sf g}_{_C}^{-1}(\xi_{_s})_{_{(1)}}$ can be described as follows:
\beqa\label{bigstratum}
{\sf g}_{_C}^{-1}(\xi_{_s})_{_{(1)}} \cap \nu_{_*}^{^{-1}}(\eta) \simeq \frac{\overbrace{\bg_{_m} \times \ldots \times \bg_{_m}}^{n}}{\Delta(\bg_{_m})}.
\eeqa
\eprop
\bpr By Proposition \ref{ss spectral curve} we have an identification 
\beqa
{\sf g}_{_C}^{-1}(\xi_{_s})_{_{(1)}} \simeq G^{^{(1)}}(1,\delta)/\Delta(\bg_{_m})
\eeqa
and by Proposition \ref{oscaporaso}, we see that 
\beqa
{P_{_{\delta,{Y_{_{\xi,s}}}}}}(n) \simeq G^{^{(1)}}(1,\delta)/\underbrace{\bg_{_m} \times \ldots \times \bg_{_m}}_{n}.
\eeqa
Thus we get the required identification \eqref{bigstratum} of the big stratum of the general Gieseker-Hitchin fibre.\epr

\brem A similar description clearly holds for the other strata as well. \erem
\brem Proposition \ref{opentoric} should be viewed in the light of the following remarks.  Let $E$ be a torsion-free $\co_{_C}$-module such that the local structure at the node on $C$ is of type ${\mathfrak m}^{^n}$. Then by \cite[Remark 5.2]{ictp},  the fibre $\nu_{_*}^{^{-1}}(E)$ can be identified with the so-called {\em wonderful compactification} of $PGL(n)$. 

Consider ${P_{_{\delta,{Y_{_{\xi,s}}}}}}$ the compactified Picard variety of the spectral curve $Y_{_{\xi,s}}$ over the closed fibre $C$.  We view ${P_{_{\delta,{Y_{_{\xi,s}}}}}}$ as a subscheme of ${\eM}_{_{C}}^{^{H}}(n,d)$ of torsion-free Hitchin pairs on $C$. Under this identification, a point $\eta \in {P_{_{\delta,{Y_{_{\xi,s}}}}}}(n)$ gives a torsion-free Hitchin pair $(E, \theta)$, such that 
the local structure at the node on $C$ is of type ${\mathfrak m}^{^n}$. 

Proposition \ref{opentoric} shows that the inclusion 
\beqa {\sf g}_{_C}^{-1}(\xi_{_s})_{_{(1)}} \cap \nu_{_*}^{^{-1}}(\eta) \subset \nu_{_*}^{^{-1}}(E)
\eeqa 
is in fact the inclusion 
\beqa
\frac{\overbrace{\bg_{_m} \times \ldots \times \bg_{_m}}^{n}}{\Delta(\bg_{_m})} \subset PGL(n)
\eeqa
i.e the standard inclusion of the maximal torus of $PGL(n)$. 
\erem

Let $\tilde{g} = n(g - 1) + 1 + deg(\cl).\frac{n(n-1)}{2}$ be the arithmetic genus of the spectral vine curve $Y_{_{\xi,s}}$ (see Remark \ref{spectralproperty}). 

By Corollary \ref{asnc},  the structure morphism Gieseker-Hitchin scheme ${\eG}_{_{S}}^{^{H}}(n,d) \to S$ is flat and the closed fibre ${\eG}_{_{C}}^{^{H}}(n,d)$ has analytic normal crossing singularities. Because we have shown the Gieseker-Hitchin morphism  is proper and because we are in characteristic zero, the general fibre also has analytic normal crossing singularities.

In summary we have proven the following main theorem.
\bth\label{relative abelianization}(\it Quasi-abelianization) Consider the restriction of the proper birational morphism $\nu_{_*}:{\sf g}_{_C}^{-1}(\xi_{_s}) \to {P_{_{\delta,{Y_{_{\xi,s}}}}}}$.   
\begin{enumerate}
\item Let $\eta \in {P_{_{\delta,{Y_{_{\xi,s}}}}}}(j)$. The fibre $\nu_{_*}^{^{-1}}(\eta)$ can be identified with the projective toric variety $\overline {T_{_j}}$ which is the closure of the maximal torus $T_{_j} \subset PGL(j)$ in the wonderful compactification $\overline {PGL(j)}$. This is in fact the toric variety associated to the Weyl chamber of $PGL(j)$ (\cite{procesi}). 
\item The scheme ${\sf g}_{_S}^{-1}(\xi)$ provides a relative compactification of the Jacobian of smooth (spectral) curves of genus $\tilde{g}$ which has a divisor ${\sf g}_{_C}^{-1}(\xi_{_s})$ with analytic normal crossing singularities. 
\end{enumerate}
\eeth
\brem\label{morecomparison} The last statement poses the interesting general problem of giving a modular construction of a compactified Picard variety for stable curves which has analytic normal crossing singularities. For the case of a vine curve with $n$-nodes $n \geq 2$, one needs to consider quasi-Gieseker line bundles on curves $\{Y^{^{(\ell)}}\}_{_{\ell = 1}}^{^n}$, unlike the Caporaso compactification which requires only quasi-Gieseker line bundles on the {\it ladder} curve alone. \erem

\section{The reducible curve case}
The above theory goes through in the case when the closed fibre of $X \to S$ is a reducible curve  $C = C_{_1} \cup C_{_2}$  with a single node at a point $p \in C_{_1}\cap C_{_2}$,
where $C_{_1}$ and $C_{_2}$ are two smooth curves over an algebraically closed field $k$ of genus $g_1$ and $g_2$ respectively. 
 
If $C$ is irreducible, and $\cl$ be an invertible sheaf over $C$ then $\cl$ is obtained by giving line bundles $\cl_i$ on $C_i$ together with a gluing isomorphism $\ell:\cl_{1,p} \simeq \cl_{2,p}$.
 
A {\em polarization} on the reducible nodal curve $C$ can be thought of as giving a pair $a = (a_1,a_2)$  with $a_i > 0$ positive rational numbers  with $a_1+a_2 =1$. Let $\cl$ be an ample invertible sheaf on $C$; this in turn gives a pair of ample invertible sheaves $\cl_i$ on $C_i$. Equivalently, we say that $\cl$ gives a  {\em polarization} on $C$  if in terms of $\cl_i$ on $C_i$, one has $\frac{deg(\cl_1)}{deg(\cl_2)} = \frac{a_1}{a_2}$.
To ensure that under the assumption $gcd(n,d) = 1$ we have the condition {\em semistable = stable} for the Hitchin pairs, we need to impose a genericity condition, namely we assume that $a_1.\chi \notin \bz$. Under these hypotheses, we will be dealing only with stable objects in this paper. Note that if the curve $C$ is irreducible, only the condition $gcd(n,d) = 1$ would do since no polarization figures in the definition of stability. We will make these assumptions in this section.

The sheaf $E$ is of {\em rank} $(n_1,n_2)$, if $rank(E_i)  = n_i$, where $E_i:= E|_{{C_i}}$. Say $E$ is of rank $n$ if $n = n_1 = n_2$. Note that for a torsion-free $\co_{_C}$--module, at least one of the $n_i \neq 0$.

For a torsion-free sheaf $E$ on $C$ and the polarization $a$, define the {\it $a$--rank} and {\it $a$--slope} of $E$ as follows:
\beqa
rk_{_a}(E):= a_1.rk(E_i) + a_2.rk(E_2)
\\
\mu_{_a}(E):= \frac{\chi(E)}{rk_{_a}(E)}, \text{if}~ rk_{_a}(E) \neq 0
\eeqa

Since $dim(C)= 1$ for us, we see immediately that
\beqa\label{pandmu}
\frac{{\tt p}(E,m)}{rk_a(E)} = m \cdot deg(\cl) + \mu_a(E).
\eeqa
where ${\tt p}(E,m) := \chi(E \otimes \cl^{^m})$.

Almost all of the general theory developed above works without change for this case also and the proofs are really no different. The only new feature which emerges is that the choice of a {\em polarization} is needed to define the notion of stability of Hitchin pairs as we saw in Section 2. This naturally leads to a suitable notion of stability of the pure sheaves on the $Z = {\bp}(\cl^* \oplus \co_{_C})$. The interesting new feature is that in the description of the Hitchin fibre the choice of polarization enters the definition of the compactified Picard variety when seen from the stand-point of Oda-Seshadri. The generic spectral curve is again a {\em vine curve} with $n$-nodes and two irreducible components.

As one knows, Caporaso's construction of the compactification of the Picard variety  does not need any polarization; this can be explained by showing that the Oda-Seshadri compactification for a generic polarization is isomorphic to Caporaso's  compactification (see \cite{alex}). We summarize the results in the reducible curve case in the following theorem whose proof we omit since it is similar to the proof of previous theorem:

\bth The Gieseker-Hitchin moduli space ${\eG}_{_{S}}^{^{H}}(n,d)$ is well-defined quasi-projective scheme, flat over $S$ with generic fibre isomorphic to the classical Hitchin space. We have a Hitchin map ${\sf g}_{_S}:{\eG}_{_{S}}^{^{H}}(n,d) \to {\ca}_{_S}$ which is {\em proper}. There is a proper birational morphism from the  Hitchin fibre over a general section $\sigma:S \to {\ca}_{_S}$ to the compactified relative Picard \text{$S$-scheme} ${P_{_{\delta,Y_{_\sigma}}}}$ of the spectral surface $Y_{_\sigma}$.  \eeth
\brem The spectral surface $Y_{_\sigma}$ is a fibered surface over $S$ with smooth generic fibre and the closed fibre is a reducible vine curve with two components and $n$-nodes. \erem 
\brem We remark that in the special case when $n=2$ and $d=1$ and $C$ is {\em reducible} with two components, in contrast to the general phenomenon discussed in the paper, it so happens that there is an isomorphism 
\beqa\label{specialiso}
{\eG}_{_{S}}^{^{H}}(2,1) \simeq {\eM}_{_{S}}^{^{H}}(2,1)
\eeqa 
in this special situation, which is somewhat misleading since in this case the Hitchin map is well-defined even on the moduli of torsion-free Hitchin pairs. The choice of a polarization allows the discarding of torsion-free sheaves which locally look like $\mathfrak m \oplus \mathfrak m$ and this is the reason for the isomorphism \eqref{specialiso} in the rank $2$ case. Details of this will appear in the second author's doctoral thesis (\cite{barik}).\erem


\begin{thebibliography}{9999}

\bibitem{alex} V. Alexeev,  Compactified Jacobians and Torelli map, {\it Publ. Res. Inst. Math. Sci.} {\bf 40} (2004), 1241-1265.

\bibitem{barik} P. Barik, Doctoral Thesis ``Higgs sheaves on a singular curve", Chennai Mathematical Institute, 2013.
\bibitem{bnr} A. Beauville, M.S. Narasimhan, S. Ramanan, Spectral curves and the generalised theta divisor, {\it J. Reine Angew. Math.} {\bf 398} (1989), 169-179.

\bibitem{bisram} Indranil Biswas, S. Ramanan, An infinitesimal study of the moduli of Hitchin pairs, {\it J. London Math.Soc.} {\bf 49} (1994), 219-231.

\bibitem{caporaso} L. Caporaso, A compactification of the universal Picard variety over the moduli space of stable curves, {\it Journal of the American  Society}, {\bf 7}, Number 3 (1994). 
\bibitem{conrad} B. Conrad, Deligne's notes on Nagata compactifications, {\it J. Ramanujan Math. Soc.}, {\bf 22} (2007), 205-257.

\bibitem{gies} D. Gieseker, A Degeneration of the Moduli Space of Stable Bundles, {\it J. Differential Geometry}, {\bf 19} (1984), 173-206.


\bibitem{ega} A. Grothendieck, J. Dieudonne, Elements de Geom\'etrie Alg\'ebrique-I, Grundlehren, Springer, 1971.
\bibitem{hartshorne} R. Hartshorne, Algebraic Geometry, GTM 52, Springer.

\bibitem{hitchin} N.J. Hitchin, The self-duality equations on a Riemann surface, {\it Proc. London Math. Soc.} , {\bf 55} (1987),  59-126. 

\bibitem{hitchin1} N.J. Hitchin, Stable bundles and integrable systems, {\it Duke Math. Journal}, {\bf 54} (1987), 91-114. 

\bibitem{kausz} I. Kausz, A Gieseker type degeneration of moduli stacks of vector bundles on curves,
{\it Transactions of the American Mathematical Society},
{\bf 357}, Number 12 (2004), 4897-4955.

\bibitem{lipman} J. Lipman, Rational singularities, with applications to algebraic surfaces and unique factorization. {\it
Pub. Math. I.H.E.S.}, {\bf 36} (1969), 195-279.


\bibitem{mumford} D. Mumford, Abelian Varieties, Oxford University Press, 2010.



\bibitem{ns2} D.S. Nagaraj and C.S. Seshadri,  Degenerations of the moduli spaces of vector bundles on curves. II (Generalized Gieseker moduli spaces). {\it Proc. Indian Acad. Sci. Math. Sci.} {\bf 109} no. 2 (1999), 165-201.






\bibitem{nitsure} N. Nitsure, Moduli space of semistable pairs on a curve, {\it Proc. London Math. Soc.} (3) {\bf 62} (1991), 275-300.

\bibitem{os} T. Oda and C.S. Seshadri, Compactifications of the generalized Jacobian variety, {\it Trans. Amer. Math. Soc.} {\bf 253} (1979), 1-90. 

\bibitem{pandharipande} R. Pandharipande, A compactification over $M_g$ of the universal moduli space of slope-semistable vector bundles. {\it J. Amer. Math. Soc.}, {\bf 9}(2)(1996),425-471.

\bibitem{procesi} Claudio Procesi, The toric variety associated to Weyl chambers, Mots, 153-161, Lang. Raison. Calc., Herms, Paris, 1990. 

\bibitem{schmitt} A. Schmitt, The Hilbert compactification of the universal moduli space of semistable vector bundles over smooth curves, {\it J. Differential Geometry}, {\bf 66} (2004), 169-209. 

\bibitem{ast} C. S. Seshadri, Fibr\'es vectoriels sur les courbes alg\'ebriques, {\it Asterisque} {\bf 96} (1982), 1-209.
\bibitem{ictp} C. S. Seshadri, Degenerations of the moduli space of vector bundles, ICTP, 1999.

\bibitem{sim} C. Simpson, Hitchin pairs and local systems, {\it Publ. Math. I.H.E.S.} {\bf 75} (1992), 5-95.

\bibitem{sim1} C. Simpson, Moduli of representations of the
fundamental group of a smooth projective variety-I, {\it
Pub. Math. I.H.E.S.} {\bf 79} (1994)  47-129.

\bibitem{sim2} C. Simpson, Moduli of representations of the
fundamental group of a smooth projective variety-II, {\it
Pub. Math. I.H.E.S.} {\bf 80} (1995), 5-79.

\bibitem{teod} Titus Teodorescu, Principal bundles over chains or cycles of rational curves, {\it Michigan Math. J.}, {\bf 50}, Issue 1 (2002), 173-186.

\bibitem{sun} X.T. Sun: Factorization of generalized theta functions in the reducible case, {\it Ark. Mat.}, {\bf 41} (2003), 165-202.

\end{thebibliography}
\end{document}